\documentclass[onecolumn]{elsarticle}
\makeatletter
\def\ps@pprintTitle{%
	\let\@oddhead\@empty
	\let\@evenhead\@empty
	\def\@oddfoot{\centerline{\thepage}}%
	\let\@evenfoot\@oddfoot}
\makeatother
\usepackage[english]{babel}
\usepackage{amsfonts}
\usepackage{amsmath}
\usepackage{amscd}
\usepackage{dirtytalk}
\usepackage[active]{srcltx} % SRC Specials: DVI [Inverse] Search
\usepackage{latexsym}
\usepackage{amssymb}
\usepackage[latin1]{inputenc}
\usepackage{enumerate}

\usepackage{pdfsync}
\usepackage{amsthm}
\usepackage{amsmath}
\usepackage{mathtools}
\usepackage{amscd}
\usepackage{latexsym}
\usepackage{amssymb}
\usepackage{enumitem}
\newtheorem{thm}{Theorem}[section]
\newtheorem{prop}[thm]{Proposition}
\newtheorem{cor}[thm]{Corollary}
\newtheorem*{cor*}{Corollary}
\newtheorem{lema}[thm]{Lemma}
\newtheorem*{lema*}{Lemma}

\numberwithin{equation}{section}
\theoremstyle{definition}
\newtheorem*{Def}{Definition}

\newenvironment{dem}{\vspace{1ex}\noindent{\it Proof.}\hspace{0.5em}}
{\hfill\qed\vspace{1ex}}

\newtheorem*{obs}{Remark}

\newtheorem*{obs*}{Remark}
\newtheorem*{thm*}{Theorem}
\newtheorem*{prop*}{Proposition}

\newtheoremstyle{dotless}{}{}{}{}{}{}{ }{}
\theoremstyle{dotless}

\newcommand{\matriz}[4]{\displaystyle\
    \left(
       \begin{array}{cc}
        {#1}&{#2}\\
        {#3}&{#4}
       \end{array}
     \right)}

 \newcommand{\mm}[3]{\displaystyle\
 	\begin{array}{c}
 		{#1}\\
 		{#2} \\
 		{#3}
 	\end{array}
 }

 \newcommand{\ms}[2]{\displaystyle\
	\begin{array}{c}
		{#1}\\
		{#2}
	\end{array}
}

\newcommand{\PI}[2]{\left\langle \,#1 , #2\, \right\rangle}

\newcommand{\K}[2]{[ \,#1 , #2\, ]}

\newcommand{\ra}{\rightarrow}

\newcommand{\WS}{W_{/ [\St]}}

\newcommand{\vareps}{\varepsilon}

\newcommand{\St}{\mathcal{S}}

\newcommand{\HH}{\mathcal{H}}

\newcommand{\M}{\mathcal{M}}
\newcommand{\N}{\mathcal{N}}

\newcommand{\Q}{\mathcal{Q}}
\newcommand{\T}{\mathcal{T}}

\newcommand{\KK}{\mathcal{K}}
\newcommand{\mc}[1]{\mathcal{#1}}

\newcommand{\ol}{\overline}

\newcommand{\perpi}{[\perp]}

\begin{document}

\begin{frontmatter}
	
	\title{Schur complements of selfadjoint Krein space operators}
	%\title{Approximate solutions of linear operators equations in weighted Schatten $p$-norms}
	
	\author[FI,IAM]{Maximiliano Contino\corref{ca}}
	\ead{mcontino@fi.uba.ar}
	
	\author[FI,IAM]{Alejandra Maestripieri}
	\ead{amaestri@fi.uba.ar}
	
		\author[IAM,IV,UNGS]{Stefania Marcantognini}
	\ead{smarcantognini@ungs.edu.ar}

	\cortext[ca]{Corresponding author}
	\address[FI]{%
		Facultad de Ingenier\'{\i}a, Universidad de Buenos Aires\\
		Paseo Col\'on 850 \\
		(1063) Buenos Aires,
		Argentina 
	}
	
	\address[IAM]{%
		Instituto Argentino de Matem\'atica ``Alberto P. Calder\'on'' \\ CONICET\\
		Saavedra 15, Piso 3\\
		(1083) Buenos Aires, 
		Argentina }
	
	\address[IV]{%
		Departamento de Matem\'atica -- Instituto Venezolano de Investigaciones Cient\'ificas \\ Km 11 Carretera Panamericana \\ Caracas, Venezuela
	}

	\address[UNGS]{%
 Universidad Nacional de General Sarmiento -- Instituto de Ciencias \\ Juan Mar\'ia Gutierrez \\ (1613) Los Polvorines, Pcia. de Buenos Aires, Argentina
}
	
	\begin{abstract}
	Given a bounded selfadjoint operator $W$ on a Krein space $\HH$ and a closed subspace $\St$ of $\HH$,  the Schur complement of $W$ to $\St$ is defined under the hypothesis of weak complementability. A variational characterization of the Schur complement is given and the set of selfadjoint operators $W$ admitting a Schur complement with these variational properties is shown to coincide with the set of $\St$-weakly complementable selfadjoint operators.
	\end{abstract}
	
	
\end{frontmatter}

\section{Introduction}
The notion of Schur complement (or shorted operator) of $B$ to $\St$ for a positive operator $B$ on a Hilbert space $\HH$ and $\St \subseteq \HH$ a closed subspace, was introduced by M.G.~Krein \cite{Krein}. When $\leq_{\HH}$ is the usual order in $L(\HH)$, he proved that the set $\{ X \in L(\HH): \ 0\leq_{\HH} X\leq_{\HH} B \mbox{ and } R(X)\subseteq \St^{\perp}\}$ has a maximum element, which he defined as the {{Schur complement}} $B_{/ \St}$ of $B$ to $\St.$ This notion was later rediscovered by Anderson and Trapp \cite{Shorted2}. If $B$ is represented as the $2\times 2$ block matrix $\begin{pmatrix} a & b\\ b^* & c \end{pmatrix}$ with respect to the decomposition of $\HH = \St \oplus \St^{\perp},$ they established the formula 
$$B_{/ \St}= \begin{pmatrix} 0 & 0\\ 0& c - y^*y\end{pmatrix}$$ where $y$ is the unique solution of the equation $b = a^{1/2} x$ such that the range inclusion $R(y) \subseteq \ol{R(a)}$ holds. The solution always exists because $B$ is positive, in which case $a$ is also positive and the range inclusion $R(b) \subseteq R(a^{1/2})$ holds.

In \cite{AntCorSto06} Antezana et al., extended the Schur complement to any bounded operator $B$ satisfying a weak complementability condition with respect to a given pair of closed subspaces $\St$ and $\T,$ by giving an Anderson-Trapp type formula. In particular, if $B$ is a bounded selfadjoint operator, $\St=\T$ and $B=\begin{pmatrix} a & b\\ b^* & c \end{pmatrix},$ this condition reads $R(b) \subseteq R(\vert a \vert^{1/2}),$ which as noted, is automatic for positive operators. 
Later, Massey and Stojanoff \cite{Massey} studied many properties of the Schur complement of an $\St$-weakly complementable selfadjoint operator $B$ when $\St$ is $B$-positive.

In this paper we show that if $B$ is a bounded selfadjoint operator which is $\St$-weakly complementable then $B_{/ \St}$ can be characterized as the solution of a $min-max$ problem, extending the original approach of Krein. But, more importantly, the converse is true, in the sense that if the solution of this $min-max$ problem exists then $B$ has to be $\St$-weakly complementable. In other words, the $\St$-weakly complementable operators are exactly those selfadjoint operators admitting a Schur complement that satisfies these variational properties. 

A closed-form expression for the Schur complement $B_{/ \St}$ of $B$ to $\St$ is also established, in terms of a family of densely defined projections with prescribed nullspace $\St^{\perp}$ (Theorem \ref{ThmWC5}). This formula is new even in the case of positive $B.$   

We then turn to the consideration of a bounded selfadjoint operator $W$ on a Krein space $(\HH, \K{ \ }{ \ }).$ For a fixed signature operator $J$ on $\HH$, $JW$ is selfadjoint in the Hilbert space inner product 
$\PI{ \ }{ \ }$ associated with $J$.  If $\St$ is a given closed subspace of $\HH$, $JW$ is assumed to be $\St$-weakly complementable and $J_\alpha$ is any other signature operator on $\HH$ then two key results are established: $J_\alpha W$ is $\St$-weakly complementable (Theorem \ref{WC}) and $J~(JW)_{/ \St} = J_\alpha ~ (J_\alpha W)_{/ \St}$ (Theorem \ref{thethm}). 

Based on these results we extend the notions of $\St$-weak complementability and Schur complement to the Krein space setting. A bounded selfadjoint operator $W$ on a Krein space $\HH$ is $\St$-weakly complementable if, for some (and, hence, any) signature operator $J$, $JW$ is $\St$-weakly complementable in the corresponding Hilbert space. If this is the case then the Schur complement of $W$ to $\St$ is  $W_{/ [\St]} := J~(JW)_{/ \St}$.

In this fashion we obtain a simple way of computing the Schur complement of $\St$-weakly complementable selfadjoint operators in Krein spaces. This definition allows us to ``translate'' the properties obtained in Hilbert spaces to the Krein space setting in a straightforward way. 

If $\St$ is a regular subspace of $\HH$ (meaning that $\HH=\St \ [\dotplus] \ \St^{\perpi})$ then it is possible to give a characterization of the $\St$-weak complementability of $W$ in terms of the entries of the first row of the $2\times 2$ block matrix representation of $W$ with respect to $\St \ [\dotplus] \ \St^{\perpi}.$ Indeed if 
$W=\left(w_{ij}\right)_{i,j =1,2}$ and $w_{11} = dd^\#$ is a Bogn\'ar-Kr\'amli factorization of $w_{11}$ obtained as in \cite[Theorem 1.1]{DR1}, then $W$ is $\St$-weakly complementable if, and only if, $R(w_{12}) \subseteq R(d)$. In this case, $W_{/ [\St]}= \begin{pmatrix} 0 & 0\\ 0& w_{22} - y^{\#}y\end{pmatrix}$ with $y$ the only solution of the equation $w_{12} = d x$. The result may be viewed as a neat Krein space counterpart of the Hilbert space results in \cite{AntCorSto06}.

Based on a formula given by Pekarev \cite{Pekarev}, Maestripieri and Mart\'inez Per\'ia \cite{SchurKrein} extended the concept of the Schur complement to bounded selfadjoint operators in Krein spaces with the so-called ``unique factorization property''. Another approach was given by Mary \cite{XavierMary}. He defined the Schur complement of a bounded operator $W=\left(w_{ij}\right)_{i,j =1,2}$ when the range $R(w_{11})$ and the nullspace $N(w_{11})$ of $w_{11}$ are regular subspaces. The approach we adopt has greater scope and is less restrictive.

The paper has three additional sections. Section 2 is a brief expository introduction to Krein spaces and operators on them and serves to fix the notation and give some results that are needed in the following sections. Section 3 is entirely devoted to the study of complementability and the Schur complement of a selfadjoint operator on a Hilbert space.  In Section 4 we present our main results concerning  the Schur complement of a Krein space operator. This section includes three subsections: the first deals with the notion of weak complementability on Krein spaces; the second presents an application inspired on some completion problems previously considered in Hilbert and Krein spaces by  Baidiuk and Hassi in \cite{Baidiuk} and \cite{Baidiuk2}; in the last subsection our notion of Schur complement in the Krein space setting is compared to those in \cite{SchurKrein} and \cite{XavierMary}.

\section{Preliminaries}

We assume that all Hilbert spaces are complex and separable. If $\HH$ and $\KK$ are Hilbert spaces, $L(\HH, \KK)$ stands for the space of all the bounded linear operators from $\HH$ to $\KK.$ When $\HH = \KK$ we write, for short, $L(\HH).$ 
The domain, range and nullspace of any given $A \in L(\HH, \KK)$ are denoted by $Dom(A),$ $R(A)$ and $N(A)$, respectively. Given a subset $\T \subseteq \KK,$ the preimage of $\T$ under $A$ is denoted by $A^{-1}(\T)$ so $A^{-1}(\T)=\{ x \in \HH: \ Ax \in \T \}.$

The direct sum of two closed subspaces $\M$ and $\N$ of $\HH$ is represented by $\M \dot{+} \N.$ 
If $\HH$ is decomposed as $\HH=\M \dot{+} \N,$ the projection onto $\M$ with nullspace $\N$ is denoted by $P_{\M {\mathbin{\!/\mkern-3mu/\!}} \N}$ and abbreviated $P_{\M}$ when $\N = \M^{\perp}.$ In general, $\Q$ indicates  the subset of all the oblique projections in  $L(\HH),$ namely, $\Q:=\{Q \in L(\HH): Q^{2}=Q\}.$ 

$L(\HH)^{s}$ stands for the set of selfadjoint operators in $ L(\HH).$ Denote by $GL(\HH)$ the group of invertible operators in $L(\HH),$ $L(\HH)^+$ the cone of positive semidefinite operators in  $L(\HH)$ and $GL(\HH)^+:=GL(\HH) \cap L(\HH)^+.$
Given two operators $S, T \in L(\HH),$ the notation  $T \leq_{\HH} S$ signifies that $S-T \in L(\HH)^+.$ 
Given any $T \in L(\HH),$ $\vert T \vert := (T^*T)^{1/2}$ is the modulus of $T$ and $T=U\vert T\vert$ is the polar decomposition of  $T,$ with $U$ the partial isometry such that $N(U)=N(T).$

\vspace{0.35cm}
The following  is a well-known result about range inclusion and factorizations of operators.
\begin{lema} [Douglas' Lemma \cite{Douglas}] \label{teo Douglas} Let $Y \in L(\KK_1, \HH)$ and $Z \in L(\KK_2, \HH)$. Then
	$R(Z)\subseteq R(Y)$ if and only if there exists $D\in L(\KK_2, \KK_1)$ such that $Z=YD.$
	
	Amongst the solutions of the equation $Z=YX,$ there exists a unique operator $D_0 \in L(\HH)$ such that $N(Z)=N(D_0)$ and $R(D_0) \subseteq \ol{R(Y^*)}.$ 
\end{lema}
The operator $D_0$  is called the \emph{reduced solution} of $Z=YX.$

\vspace{0.35cm}
Given $B \in L(\HH)$ selfadjoint and $\St$ a closed subspace of $\HH,$ we say that $\St$ is $B$-\emph{positive} if 
$\PI{Bs}{s} > 0 \mbox{ for every } s \in \St, \ s\not =0.$ $B$-\emph{nonnegative}, $B$-\emph{neutral}, $B$-\emph{negative} and $B$-\emph{nonpositive} subspaces are defined analogously. If $\St$ and $\T$ are two closed subspaces of $\HH,$ the notation $\St \ \oplus_{B} \  \T$ is used to indicate the orthogonal direct sum of $\St$ and $\T$ when, in addition, $\PI{Bs}{t}=0 \mbox{ for every } s \in \St \mbox{ and } t \in \T.$ 

The following is a consequence of the spectral theorem for Hilbert space selfadjoint operators. 

\begin{lema} \label{lemaWdecom} Let $B \in L(\HH)^s$ and let $\St$ be a closed subspace of $\HH.$  Then $\St$ can be represented as 
	\begin{equation} \label{WdecompS}
	\St = \St_{+} \ \oplus_{B} \ \St_{-}
	\end{equation}
	where $\St_{+}$ and $\St_-$ are closed, $\St_+$ is $B$-nonnegative, $\St_-$ is $B$-nonpositive.  
\end{lema}

Let 
\begin{equation} \label{matrixWa}
B=\begin{bmatrix}
a  & b \\ 
b^* & c \\
\end{bmatrix} \ms{\St}{\St^{\perp}}
\end{equation}
be the matrix decomposition of $B$ induced by $\St$ and consider $	\St = \St_{+} \ \oplus_{B} \ \St_{-}$ as in \eqref{WdecompS}. Then, the matrix representations of $a, \vert a \vert, \vert a \vert^{1/2} \in L(\St)^s$ induced by $\St_{+}$ are:
$a= \begin{bmatrix}
a_+ & 0 \\ 
0 & -a_-\\
\end{bmatrix},$ $\vert a \vert= \begin{bmatrix}
a_+ & 0 \\ 
0 & a_-\\
\end{bmatrix},$  $\vert a \vert^{1/2}= \begin{bmatrix}
a_+^{1/2} & 0 \\ 
0 & a_-^{1/2}\\
\end{bmatrix},$ respectively.

Let us write $b:=\begin{bmatrix}
b_+ \\ 
b_-\\
\end{bmatrix}: \St^{\perp} \ra \begin{bmatrix}
\St_+ \\ 
\St_-\\
\end{bmatrix},$ where $b_{\pm}=P_{\St_{\pm}}b.$ Then $B$ can be written as
\begin{equation} \label{BS+}
B= \begin{bmatrix}
a_+ & 0 & b_+ \\ 
0 & -a_- & b_-\\
b_+^* & b_-^* & c \\ 
\end{bmatrix} \mm{\St_+}{\St_-}{\St^{\perp}}.
\end{equation}

\bigskip
In \cite[Theorem 3]{Shorted2}, the Schur complement $B_{/ \St}$ of an operator $B \in L(\HH)^+$ was characterized in the following fashion: if the matrix representation of $B$ is given by \eqref{matrixWa}, then 
$$B_{/ \St}= \begin{bmatrix}
0 & 0 \\ 
0 & c-f^*f \\
\end{bmatrix},$$ where $f$ is the reduced solution of $a^{1/2} x=b$ (which always exists for positive operators). The next lemma characterizes the positive operators in terms of its matrix decomposition. It follows easily from the fact that $B - B_{/ \St}  \geq_{\HH} 0.$ 

\begin{lema} \label{LemmaPositive} Let $\St \subseteq \HH$ be a closed subspace and $B \in L(\HH)^s$  with matrix decomposition induced by $\St$ as in \eqref{matrixWa}. Then $B \in L(\HH)^+$ if and only if 
	$$a \geq_{\HH}0, \, \,  b=b^*, \, \,  R(b) \subseteq R(a^{1/2}), \mbox{ and } c=f^*f+t,$$ with $f$ the reduced solution of the equation $b=a^{1/2}x$ and $t \geq_{\HH}0.$
\end{lema}

\subsection{Krein Spaces}

Although familiarity with operator theory on Krein spaces is presumed, we include some basic notions. Standard references on Krein spaces and operators on them are \cite{AndoLibro}, \cite{Azizov} and \cite{Bognar}. We also refer to \cite{DR} and \cite{DR1} as authoritative accounts of the subject.

Consider a linear space $\HH$ with an indefinite metric; i.e., a sesquilinear Hermitian form $\K{ \ }{ \ }.$  A vector $x \in \HH$ is said to be {\emph{positive}} if $\K{x}{x} > 0.$ A subspace $\St$ of $\HH$ is {\emph{positive}} if every $x \in \St,$ $x \not =0,$ is a positive vector. {\emph{Negative, nonnegative, nonpositive}} and {\emph{neutral}} vectors and subspaces are defined likewise. 

We say that two closed subspaces $\M$ and $\N$ are {\emph{orthogonal}}, and  write $\M  \ \perpi \ \N,$ if $\K{m}{n}=0 \mbox{ for every } m \in \M \mbox{ and } n \in \N.$ Denote the orthogonal direct sum of two closed subspaces $\M$ and $\N$ by $\M \ [\dotplus] \ \N.$  

Given any subspace $\St$ of $\HH,$ the {\emph{orthogonal companion}} of $\St$ in $\HH$ is defined as 
$$\St^{\perpi}:=\{ x \in \HH: \K{x}{s}=0 \mbox{ for every } s \in \St \}.$$

An indefinite metric space $(\HH, \K{ \ }{ \ })$ is a {\emph{Krein space}} if it admits a decomposition as an orthogonal direct sum 
\begin{equation} \label{fundamentaldecom}
\HH=\HH_+ \ [\dotplus] \ \HH_-,
\end{equation} 
where $(\HH_+, \K{ \ }{ \ })$ and $(\HH_-, -\K{ \ }{ \ })$ are Hilbert spaces. Any decomposition with these properties is called a {\emph{fundamental decomposition}} of  $\HH.$ 
%Notice that any Hilbert space is itself a Krein space.

Given a Krein space $(\HH, \K{ \ }{ \ })$ with a fundamental decomposition $\HH=\HH_+ \ [\dotplus] \ \HH_-,$ the (orthogonal) direct sum of the Hilbert spaces $(\HH_+, \K{ \ }{ \ })$ and $(\HH_-, -\K{ \ }{ \ })$ is a Hilbert space, $(\HH, \PI{ \ }{ \ }).$ Notice that the inner product $\PI{ \ }{ \ }$ and the corresponding quadratic norm $\Vert \ \Vert$ depend on the fundamental decomposition. 

Every fundamental decomposition of $\HH$ has an associated {\emph{signature operator}}: $J:=P_{+} - P_{-}$ with $P_{\pm}:=P_{\HH_{\pm}  {\mathbin{\!/\mkern-3mu/\!}} \HH_{\mp}}.$ The indefinite metric and the inner product corresponding to a fundamental decomposition of $\HH$ with signature operator $J$ are related to each other by 
$$\PI{x}{y}=\K{Jx}{y} \quad (x, y \in \HH).$$ 

If $\HH$ is a Krein space, $L(\HH)$ stands for the vector space of  all the linear operators on $\HH$ which are bounded in an associated Hilbert space $(\HH, \PI{ \ }{ \ }).$ Since the norms generated by different fundamental decompositions of a Krein space $\HH$ are equivalent (see, for instance, \cite[Theorem 7.19]{Azizov}), $L(\HH)$ does not depend on the chosen underlying Hilbert space, all of which are equivalent. 

Given $T \in L(\HH),$ $T^{\#}$ is the unique operator satisfying $$\K{Tx}{y}=\K{x}{T^{\#}y} \mbox{ for every } x, y \in \HH.$$
$L(\HH)^{[s]}$ denotes the set of the operators $T \in L(\HH)$ such that $T=T^{\#}.$ The selfadjoint operator $T\in L(\HH)$ is {\emph{positive}} if $\K{Tx}{x} \geq 0 \mbox{ for every } x \in \HH.$ The notation  $S \leq T$ signifies that $T-S$ is positive.

A (closed) subspace $\St$ of a Krein space $\HH$ is {\emph{regular}} if it is itself a Krein space in the indefinite metric of $\HH.$ A subspace $\St$ is regular if and only if $\HH=\St \ [\dotplus] \ \St^{\perpi}$ or, equivalently, if it is the range of a selfadjoint projection, i.e., there exists $Q \in \Q$ such that $Q=Q^{\#}$ and $R(Q)=\St$ (see \cite[Proposition 1.4.19]{Azizov}). Clearly, $\St$ is regular if and only if $\St^{\perpi}$ is regular.

Suppose that $\St$ is a regular subspace with fundamental decomposition $\St=\St_+ \ [\dotplus] \ \St_-.$ Then, by \cite[Theorem 1.6]{DR1}, there exists a fundamental decomposition of $\HH= \HH_+ \  [\dotplus] \ \HH_-$ such that $\St_{\pm} \subseteq \HH_{\pm}.$ In this case, $$\HH_{\pm}=\St_{\pm} \  [\dotplus] \ \N_{\pm},$$ where $\St^{\perpi}=\N_+ \  [\dotplus]  \ \N_-$ is a fundamental decomposition of $\St^ {\perpi}.$ Now, consider $J_1=\matriz{I}{0}{0}{-I} \ms{\St_+}{\St_-}$ and $J_2=\matriz{I}{0}{0}{-I} \ms{\N_+}{\N_-},$ signature operators of $\St$ and $\St^ {\perpi},$ respectively. 
Then 
\begin{equation} \label{signature}
J=\matriz{J_1}{0}{0}{J_2}  \ms{\St}{\St^{\perpi}}
\end{equation}
is a signature operator for $\HH.$

Given $W \in L(\HH)^s$ and $\St$ a closed subspace of $\HH,$ we say that $\St$ is $W$-\emph{positive} if 
$\K{Ws}{s} > 0$ for every $s \in \St, \ s\not =0.$ $W$-\emph{nonnegative}, $W$-\emph{neutral}, $W$-\emph{negative} and $W$-\emph{nonpositive} subspaces are defined likewise. If $\St$ and $\T$ are two closed subspaces of $\HH,$ the notation $\St \ [\dotplus]_{W} \  \T$ is used to indicate the direct sum of $\St$ and $\T$ when, additionally, $\K{Ws}{t}=0 \mbox{ for every } s \in \St \mbox{ and } t \in \T.$ 

\section{Complementability and Schur complement for selfadjoint operators in Hilbert spaces} \label{sec3}

The notion of complementability of an operator $B \in L(\HH)$ with respect to two given closed subspaces $\St$ and $\T$ of $\HH$ was studied for matrices by Ando \cite{AndoSchur} and extended to operators in Hilbert spaces by Carlson and Haynsworth \cite{Carlson}. 
In \cite{AntCorSto06} Antezana et al. defined a weaker notion, that of \emph{weak complementability}, and extended the notion of the Schur complement to this context. We use these ideas when $\St=\T$ and $B \in L(\HH)^s.$ In what follows, we recall both definitions for this particular case:

\begin{Def} Let $B \in L(\HH)^s$ and let $\St \subseteq \HH$ be a closed subspace. Then $B$ is $\St$\emph{-complementable} if
	$$\HH=\St+ B^{-1}(\St^{\perp}).$$
\end{Def}

In \cite{CMSSzeged} it was shown that $B$ is $\St$-complementable if and only if there exists a $B$-selfadjoint projection onto $\St;$ i.e., the set
$$\mc{P}(B,\St):=\{ Q \in \Q: R(Q)=\St, \ BQ=Q^{*}B \}$$ is not empty.
It was also proven that, if
\begin{equation} \label{matrixW}
B=\begin{bmatrix}
a  & b \\ 
b^* & c \\
\end{bmatrix} \ms{\St}{\St^{\perp}},
\end{equation}
then $B$ is $\St$-complementable if and only if $R(b) \subseteq R(a).$

This naturally leads to the following definition.
\begin{Def} Let  $\St \subseteq \HH$ be a closed subspace, and $B \in L(\HH)^s$ with representation as in \eqref{matrixW}. Then $B$ is $\St$\emph{-weakly complementable} if
	$$R(b) \subseteq R(\vert a \vert^{1/2}).$$
\end{Def}
When $R(a)$ is closed both notions coincide and therefore the notion of weak complementability is distinct only in the infinite dimensional setting.
Every positive operator $B$ is $\St$-weakly complementable.

\begin{prop} \label{WCpositive} Let $B \in L(\HH)^s.$ Then $B$ is $\St$-weakly complementable for every closed subspace $\St \subseteq \HH$ if and only if $B$ is semidefinite. 
\end{prop}

\begin{dem} If $B$ is semidefinite then $B$ is $\St$-weakly complementable for every closed subspace $\St \subseteq \HH,$ because in this case, if $B \in L(\HH)$ is represented as in \eqref{matrixW} for any $\St$, by Lemma \ref{LemmaPositive}, $R(b) \subseteq R((\pm a)^{1/2}).$
	
	Conversely, suppose that $B$ is $\St$-weakly complementable for every closed subspace $\St \subseteq \HH$ and that $B$ is not definite. Then there exists $x_0 \in \HH \setminus \{0\}$ such that $x_0$ is $B$-neutral and $x_0 \not \in N(B).$ Let $\St= span \ \{ x_0\}$ and suppose that $B \in L(\HH)^s$ is represented as in \eqref{matrixW}. Then $\PI{By}{y}=\PI{ay}{y}=0 \mbox{ for every } y \in \St.$ Hence $a = 0$ and  $b=0=b^*,$ because $B$ is $\St$-weakly complementable. Then, $\St \subseteq N(B)$ which is a contradiction. Therefore, $B$ is semidefinite.
\end{dem}

Also, $B$ is $\St$-weakly complementable and $\St$ is $B$-nonnegative if and only if $a \in L(\St)^+$ and $R(b) \subseteq R(a^{1/2}).$ 
In fact, if $\St$ is $B$-nonnegative then, for every $s \in \St,$
$$0 \leq \PI{Bs}{s}= \PI{as}{s},$$ whence $a \in L(\St)^+$ and, since $B$ is $ \St$-weakly complementable, $R(b) \subseteq R(a^{1/2}).$ The converse is similar. 
Analogously, $B$ is $\St$-weakly complementable and $\St$ is $B$-nonpositive if and only if $-a \in L(\St)^+$ and $R(b) \subseteq R((-a)^{1/2}).$

We recall the definition of Schur complement for an $\St$-weakly complementable selfadjoint operator.

\begin{Def}   Let $\St \subseteq \HH$ be a closed subspace and let $B \in L(\HH)^s$ be $\St$-weakly complementable. When $B$ is as in \eqref{matrixW}, let $f$ be the reduced solution of $b=\vert a\vert^{1/2}x$ and  $a=u \vert a \vert$ the polar decomposition of $a.$ The \emph{Schur complement} of $B$ to $\St$ is defined as
	$$B_{ / \St}:= \begin{bmatrix}
	0 & 0 \\ 
	0 & c-f^*uf \\
	\end{bmatrix}.$$ 
	$B_{\St} := B- B_{/ \St}$ is the $\St$-\emph{compression} of $B.$
\end{Def}

If $B$ is positive, $B_{ / \St}$ coincides with the usual Schur complement of $B$ to $\St.$

\subsection{Variational characterization of the Schur complement}

For $B \in L(\HH)^s$ and $\St$ a closed subspace of $\HH,$ define
$$\M^{-}(B,\St^{\perp}):=\{ X \in L(\HH)^s, \ X \leq_{\HH} B, \ R(X) \subseteq \St^{\perp} \},$$
$$\M^{+}(B,\St^{\perp}):=\{ X \in L(\HH)^s, \ B \leq_{\HH} X, \ R(X) \subseteq \St^{\perp} \}.$$

The next proposition shows that $B$ is $\St$-weakly complementable if and only if $\M^{-}(B,\St_+^{\perp})$ and $\M^{+}(B,\St_-^{\perp})$ are non-empty, where $\St = \St_{+} \ \oplus_{B} \ \St_{-}$ is any decomposition as in \eqref{WdecompS}.

\begin{prop} \label{PropWC3} Let $B \in L(\HH)^s$ and let $\St$ be a closed subspace of $\HH.$ Suppose that $\St = \St_{+} \ \oplus_{B} \ \St_{-}$ is any decomposition as in \eqref{WdecompS}. Then, the following statements are equivalent:
	\begin{enumerate}
		\item[i)] $B$ is $\St$-weakly complementable;
		\item [ii)] there exist $B_1, B_2, B_3 \in L(\HH)^s,$ $B_2, B_3 \geq_{\HH} 0$ such that $B=B_1+B_2-B_3$ and $\St \subseteq N(B_1),$ $\St_- \subseteq N(B_2),$ $\St_+ \subseteq N(B_3);$  
		\item [iii)] the sets $\M^{-}(B,\St_{+}^{\perp})$ and $\M^{+}(B,\St_{-}^{\perp})$ are non-empty; 
		\item [iv)] $B$ is $\St_{\pm}$-weakly complementable.
	\end{enumerate}
\end{prop}

\begin{dem}  $i) \Rightarrow ii):$  Let $\St = \St_{+} \ \oplus_{B} \ \St_{-}$ be any decomposition as in \eqref{WdecompS}. Suppose that the matrix representation of $B$ induced by $\St_+$ is as in \eqref{BS+} and $B$ is $\St$-weakly complementable. Then $R(b_+) \subseteq R(a_+^{1/2})$ and $R(b_-) \subseteq R(a_-^{1/2}).$ In fact, since $R(b) \subseteq R(\vert a \vert^{1/2}),$  for every $y \in \St^{\perp},$ there exists $s \in \St$ such that $by=\vert a \vert^{1/2} s.$ Therefore, for every $y \in \St^{\perp},$ $b_{\pm}y=P_{\St_{\pm}}by = P_{\St_{\pm}} \vert a \vert^{1/2} s=a_{\pm}^{1/2}s.$
	Let $f$ be the reduced solution of $b_+ = a_+^{1/2}x$ and $g$ the reduced solution of $b_-= - a_-^{1/2}x.$ 
	Set $B_1:= \begin{bmatrix}
	0 & 0 & 0 \\ 
	0 & 0 & 0 \\
	0 & 0 & c-f^*f+g^*g \\ 
	\end{bmatrix} \mm{\St_+}{\St_-}{\St^{\perp}}, $ $B_2:= \begin{bmatrix}
	a_+ & 0 & b_+ \\ 
	0 & 0 & 0\\
	b_+^* & 0 & f^*f \\ 
	\end{bmatrix} \mm{\St_+}{\St_-}{\St^{\perp}}$ and $B_3:=\begin{bmatrix}
	0 & 0 & 0 \\ 
	0 & a_- & -b_-\\
	0 & -b_-^* & g^*g \\ 
	\end{bmatrix} \mm{\St_+}{\St_-}{\St^{\perp}}.$ Then $B=B_1+B_2-B_3,$ $\St \subseteq N(B_1),$ $S_- \subseteq N(B_2),$ $S_+\subseteq N(B_3)$ and, by Lemma \ref{LemmaPositive}, $B_2, B_3 \geq_{\HH} 0.$

	$ii) \Rightarrow iii):$ Since  $B_1+B_2=B+B_3 \geq_{\HH} B$ and $R(B_1+B_2) \subseteq \St_-^{\perp},$ $B_1+B_2 \in \M^+(B,\St_{-}^{\perp}).$
	Similarly, $B_1-B_3 \in \M^-(B,\St_{+}^{\perp}).$
	
	$iii) \Rightarrow iv):$ Suppose that $X_0 \in \M^{-}(B,\St_+^{\perp}).$ Since $R(X_0) \subseteq \St_+^{\perp},$ the matrix representation of $X_0$ induced by $\St_+$ is $
	X_0=\begin{bmatrix}
	0 & 0 \\ 
	0 & d \\ 
	\end{bmatrix} \ms{\St_+}{\St_+^{\perp}},$ for some $d \in L(\St_+^{\perp}).$ Suppose that the matrix representation of $B$ induced by $\St_+$ is $B=\begin{bmatrix}
	a' & b' \\ 
	b'^* & c'  \\
	\end{bmatrix} \ms{\St_+}{\St_+^{\perp}},$ with $a' \in L(\St_+)^+.$
	Since $$\begin{bmatrix}
	a' & b' \\ 
	b'^* & c'-d\\
	\end{bmatrix} = B-X_0 \geq_{\HH} 0,$$ by Lemma \ref{LemmaPositive}, $R(b') \subseteq R(a'^{1/2})$ and $B$ is $\St_+$-weakly complementable.
	In a similar way, $B$ is $\St_-$-weakly complementable.
	
	$iv) \Rightarrow i):$ Suppose that the matrix representation of $B$ induced by $\St_+$ is as in \eqref{BS+}, since $B$ is $\St_{\pm}$-weakly complementable, $R(b_{\pm}) \subseteq R(a_{\pm}^{1/2}).$ Thus,
	$$R(b) \subseteq R(b_+) + R(b_-) \subseteq R(a_+^{1/2}) \oplus R(a_-^{1/2})=R(\vert a \vert^{1/2}),$$ and $B$ is $\St$-weakly complementable. 
\end{dem}

The following result characterizes the weak $\St$-complementability of $B$ when $\St$ is $B$-nonnegative. A similar result holds in the $B$-nonpositive case. Several of the equivalences were also proven in \cite[Proposition 3.3]{Massey}. Nonetheless, we include the proofs for the sake of completeness.

\begin{prop} \label{PropWC1} Let $B \in L(\HH)^s$ and let $\St$ be a closed subspace of $\HH.$ Then the following statements are equivalent:
	\begin{enumerate}
		\item [i)] $\St$ is $B$-nonnegative and $B$ is $\St$-weakly complementable;
		\item [ii)] there exist $B_1, B_2 \in L(\HH)^s,$ $B_2 \geq_{\HH} 0$ such that $B=B_1+B_2$ and $\St \subseteq N(B_1);$  
		\item [iii)] the set $\M^{-}(B,\St^{\perp})$ is non-empty;
		\item [iv)] the set $\M^{-}(B,\St^{\perp})$ has a maximum element, namely,
	\end{enumerate}
	$$B_{/ \St}= max \ \M^{-}(B,\St^{\perp}).$$
\end{prop}

\begin{dem}  If $\St$ is $B$-nonnegative then, in the decomposition of $\St$ as in \eqref{WdecompS}, $\St_+=\St$ and $\St_- = \{0\}.$
	Applying Proposition \ref{PropWC3}, the equivalence $i) \Leftrightarrow ii)$ and the implication $ii) \Rightarrow iii)$ follow. 
	
	$iii) \Rightarrow i):$ Suppose that $X_0 \in \M^{-}(B,\St^{\perp}).$ For all $s \in \St,$ $\PI{Bs}{s} \geq \PI{X_0s}{s} = 0,$ because $X_0 \leq_{\HH} B$ and $R(X_0) \subseteq \St^{\perp}.$ Then $\St$ is $B$-nonnegative. In this case, applying Proposition \ref{PropWC3}, we also have that $B$ is $\St$-weakly complementable.
	
	$iii) \Leftrightarrow iv):$ Suppose that $X_0 \in \M^{-}(B,\St^{\perp}),$ then by $iii)  \Rightarrow i),$  $\St$ is $B$-nonnegative and $B$ is $\St$-weakly complementable. 
	Decompose $X_0$ as $X_0={X_0}_+ - {X_0}_-,$ with ${X_0}_{\pm} \in L(\St^{\perp})^+.$ Since ${X_0}_+ - {X_0}_- \leq_{\HH} B,$ it follows that $0 \leq_{\HH} {X_0}_+ \leq_{\HH} B+{X_0}_-.$ Thus, by \cite[Theorem 1]{Shorted2}, $$0 \leq_{\HH}{X_0}_{+}  \leq (B+{X_0}_-)_{/ \St} = B_{/ \St}+{X_0}_-,$$ where the last equality is a result of the fact that if $Z \in L(\HH)$ is selfadjoint and $R(Z) \subseteq \St^{\perp}$ then $B+Z$ is $\St$-weakly complementable and $(B+Z)_{/ \St}=B_{/ \St}+Z.$ Therefore
	$X_0 \leq_{\HH} B_{/ \St}.$ Finally, as $B_{/ \St} \in L(\HH)^s,$ $R(B_{/ \St}) \subseteq \St^{\perp}$ and, by Lemma \ref{LemmaPositive}, $B_{/ \St} \leq_{\HH}B.$ Hence  $B_{/ \St} \in \M^{-}(B,\St^{\perp}).$ Thus $B_{/ \St}= max \ \M^{-}(B,\St^{\perp}).$ 
	
	The converse is straightforward. 
	
\end{dem}

The Schur complement of $B$ to $\St$ satisfies a variational characterization as a min-max if and only if $B$ is $\St$-weakly complementable, as the following theorem shows.

\begin{thm} \label{thmminmax} Let $B \in L(\HH)^s$ and let $\St$ be a closed subspace of $\HH.$ Suppose that $\St = \St_{+} \ \oplus_{B} \ \St_{-}$ is any decomposition as in \eqref{WdecompS}.  
	Then $B$ is $\St$-weakly complementable if and only if there exist \linebreak $min \  \M^{+}(max  \ \M^{-}(B,\St_{+}^{\perp}), \St_{-}^{\perp})$ and $max \  \M^{-} (min  \ \M^{+}(B,\St_{-}^{\perp}), \St_{+}^{\perp}).$ 
	In this case,
	$$B_{/ \St}= min \  \M^{+}(max  \ \M^{-}(B,\St_{+}^{\perp}), \St_{-}^{\perp})=max \  \M^{-} (min  \ \M^{+}(B,\St_{-}^{\perp}), \St_{+}^{\perp}).$$
\end{thm}

\begin{dem} Suppose that $B$ is $\St$-weakly complementable. If the matrix representation of $B$ induced by $\St_+$ is as in \eqref{BS+}, then  $R(b_{\pm}) \subseteq R(a_{\pm}^{1/2}).$ Let $f$ be the reduced solution of $b_+ = a_+^{1/2}x$ and $g$ the reduced solution of $b_-=- a_-^{1/2}x.$ Then, by Proposition \ref{PropWC1}, 
	$$B_{/ \St_+}= max \ \M^{-}(B,\St_+^{\perp})= \begin{bmatrix}
	0 & 0 &0  \\ 
	0 & -a_- & b_-\\
	0  & b_-^* & c -f^*f\\ 
	\end{bmatrix} \mm{\St_+}{\St_-}{\St^{\perp}}.$$ Thus $B_{/ \St_+}$ is $\St_-$-weakly complementable and $\St_-$ is $B_{/ \St_+}$-nonpositive. Again by Proposition \ref{PropWC1}, 	
	$$(B_{/ \St_+})_{/ \St_-}=min \ \M^{+} \ (max  \ \M^{-}(B,\St_{+}^{\perp}), \St_{-}^{\perp}) =  \begin{bmatrix}
	0 & 0 \\
	0  & c -f^*f +g^*g\\ 
	\end{bmatrix}.$$
	In a similar way,
	$$(B_{/ \St_-})_{/ \St_+}=max \ \M^{-}  \ (min  \ \M^{+}(B,\St_{-}^{\perp}), \St_{+}^{\perp}) = \begin{bmatrix}
	0 & 0 \\
	0  & c -f^*f+g^*g\\ 
	\end{bmatrix}.$$
	
	Conversely, if there exist $min \  \M^{+}(max  \ \M^{-}(B,\St_{+}^{\perp}), \St_{-}^{\perp})$ and $max \  \M^{-} (min  \ \M^{+}(B,\St_{-}^{\perp}), \St_{+}^{\perp}),$ then the sets $\M^{-}(B,\St_{+}^{\perp})$ and $\M^{+}(B,\St_{-}^{\perp})$ are non-empty. So by Proposition \ref{PropWC3},  $B$ is $\St$-weakly complementable.
	
	In this case, notice that $$b=\begin{bmatrix}
	b_+ \\ 
	b_-\\
	\end{bmatrix}= \begin{bmatrix}
	a_+^{1/2} & 0 \\ 
	0 & a_-^{1/2}\\
	\end{bmatrix} \begin{bmatrix}
	f \\
	-g
	\end{bmatrix}=\vert a \vert^{1/2}(f-g),$$ and since 
	$$R(f-g) \subseteq \ol{R(a_+^{1/2})} \oplus \ol{R(a_-^{1/2})} =\ol{R(\vert a \vert^{1/2})},$$ $y:= f-g$ is the reduced solution of the equation $b=\vert a \vert^{1/2} x.$
	Also, if $u=\begin{bmatrix}
	I & 0 \\ 
	0 & -I\\
	\end{bmatrix} \ms{\St_+}{\St_-},$ then $a=u\vert a \vert =\vert a \vert u$ is the polar decomposition of $a.$ Therefore
	$y^*uy=\begin{bmatrix}
	f^*
	-g^*
	\end{bmatrix} \begin{bmatrix}
	I & 0 \\ 
	0 & -I\\
	\end{bmatrix} \begin{bmatrix}
	f \\
	-g
	\end{bmatrix}=f^*f-g^*g$ and
	$$B_{/ \St}= min \  \M^{+}(max  \ \M^{-}(B,\St_{+}^{\perp}), \St_{-}^{\perp})=max \  \M^{-} (min  \ \M^{+}(B,\St_{-}^{\perp}), \St_{+}^{\perp}).$$
\end{dem}

\begin{cor} \label{corminmax} Let $B \in L(\HH)^s$ and let $\St$ be a closed subspace of $\HH.$ Suppose that $\St = \St_{+} \ \oplus_{B} \ \St_{-}$ is any decomposition as in \eqref{WdecompS}. If $B$ is $\St$-weakly complementable then $$B_{/ \St}= (B_{/ \St_+})_{/ \St_-}=(B_{/ \St_-})_{/ \St_+}.$$
\end{cor}

In \cite[Theorem 5]{Shorted2}, Anderson and Trapp proved that if $B \in L(\HH)^+$ and $\St$ is a closed subspace of $\HH$ then
$$B_{/ \St}=inf \ \{ Q^{*}BQ : Q \in \Q, \ N(Q)=\St \}.$$
More generally:

\begin{thm} \label{Shortedsupinf} Let $B \in L(\HH)^s$ and let $\St$ be a closed subspace of $\HH.$ Suppose that $\St = \St_{+} \ \oplus_{B} \ \St_{-}$ is any decomposition as in \eqref{WdecompS}. 
	Then $B$ is $\St$-weakly complementable if and only if there exist $$\underset{ Q_- \in \Q,\  N(Q_-)=\St_-}{sup}  \left( \underset{ Q_+ \in \Q, \ N(Q_+)=\St_+}{inf} \  Q_-^*Q_+^*BQ_+Q_- \right) $$ and $$\underset{ Q_+ \in \Q, \ N(Q_+)=\St_+}{inf} \left( \underset{ Q_- \in \Q,\  N(Q_-)=\St_-}{sup} Q_+^*Q_-^*BQ_-Q_+\right).$$ In this case, \begin{align*}
	B_{/ \St} &=\underset{ Q_- \in \Q,\  N(Q_-)=\St_-}{sup}  \left( \underset{ Q_+ \in \Q, \ N(Q_+)=\St_+}{inf} \  Q_-^*Q_+^*BQ_+Q_- \right)\\
	&=\underset{ Q_+ \in \Q, \ N(Q_+)=\St_+}{inf} \left( \underset{ Q_- \in \Q,\  N(Q_-)=\St_-}{sup} Q_+^*Q_-^*BQ_-Q_+\right).
	\end{align*}	
\end{thm}

In order to prove Theorem \ref{Shortedsupinf}, we require the following lemmas.

\begin{lema} \label{ShortedM1} Let $B \in L(\HH)^s$ and let $\St$ be a closed subspace of $\HH.$ Suppose that $\St = \St_{+} \ \oplus_{B} \ \St_{-}$ is any decomposition as in \eqref{WdecompS}. If $B$ is $\St$-weakly complementable and $B=B_1+B_2-B_3$ is any decomposition as in Proposition \ref{PropWC3}, then 
	$$B_{ / \St}= B_1+{B_2}_{ / \St_+}-{B_3}_{ / \St_-}.$$
\end{lema}

\begin{dem} Since $B$ is $\St$-weakly complementable, by Proposition \ref{PropWC3}, $B$ is $\St_{\pm}$-weakly complementable. Then, proceeding as in the proof of Proposition \ref{PropWC1}, it can be checked that $B_1+B_2-{B_{3}}_{ / \St_-}=min \ \M^+(B,\St_{-}^{\perp}).$ 
	
	On the other hand, $$B_1+{B_2}_{ / \St_+}-{B_3}_{ / \St_-}= max \  \M^{-} (B_1+B_2-{B_{3}}_{ / \St_-}, \St_{+}^{\perp}).$$	In fact, $B_1+{B_2}_{ / \St_+}-{B_3}_{ / \St_-} \in \M^{-} (B_1+B_2-{B_{3}}_{ / \St_-}, \St_{+}^{\perp}).$ This follows since $B_1+{B_2}_{ / \St_+}-{B_3}_{ / \St_-} $ is selfadjoint, $B_1+{B_2}_{ / \St_+}-{B_3}_{ / \St_-}  \leq_{\HH} B_1+B_2-{B_{3}}_{ / \St_-}$ and, by \cite[Corollary 4]{Shorted2}, $R({B_3}_{ / \St_-}) \subseteq \ol{R(B_3)} \subseteq \St_+^{\perp}.$
	Let $Y \in \M^{-} (B_1+B_2-{B_{3}}_{ / \St_-}, \St_{+}^{\perp}),$ and decompose $Y$ as $Y={Y}_+ - {Y}_-,$ with ${Y}_{\pm} \in L(\St_+^{\perp})^+.$ Since ${Y}_+ - {Y}_- \leq_{\HH} B_1+B_2-{B_{3}}_{ / \St_-},$ 
	$$0 \leq_{\HH} Y_+ + {B_{3}}_{ / \St_-} \leq_{\HH} B_1+B_2+Y_-.$$ Now, since $R(B_1+Y_-) \subseteq  \St_{+}^{\perp},$ \cite[Theorem 1]{Shorted2} gives that $$0 \leq_{\HH} {Y}_{+} + {B_{3}}_{ / \St_-} \leq_{\HH} (B_1+B_2+Y_-)_{/ \St_+} =B_1+{B_{2}}_{ / \St_+}+Y_-;$$ i.e., 
	$Y \leq_{\HH} B_1+{B_2}_{ / \St_+}-{B_3}_{ / \St_-}.$
	Therefore, 
	$$ B_1+{B_2}_{ / \St_+}-{B_3}_{ / \St_-}=max \  \M^{-} (min  \ \M^{+}(B,\St_{-}^{\perp}), \St_{+}^{\perp})=B_{ / \St}.$$
\end{dem}

\begin{lema} \label{lemasupinf} Let $B \in L(\HH)^+$  and let $\St$ be a closed subspace of $\HH$ decomposed as in \eqref{WdecompS}. For any $Q_- \in \Q,\  N(Q_-)=\St_-,$ 
	$$\underset{ Q_+ \in \Q, \ N(Q_+)=\St_+}{inf} \  Q_-^*Q_+^*BQ_+Q_-=Q_-^*B_{/ \St_+}Q_-.$$
\end{lema}

\begin{dem} Let $Q_- \in \Q,\  N(Q_-)=\St_-.$ By \cite[Theorem 5]{Shorted2}, $B_{/ \St_+} \leq Q_+^*BQ_+,$ for every $Q_+ \in \Q, \ N(Q_+)=\St_+.$ Therefore $Q_-^*B_{/ \St_+}Q_-$ is a lower bound of the set $\{Q_-^*Q_+^*BQ_+Q_- :Q_+ \in \Q, \ N(Q_+)=\St_+\}.$
	
	If $B$ is invertible then, by \cite[Section 4]{CMSSzeged}, $B$ is $\St_+$-complementable. So, by \cite[Proposition 4.2]{CMSSzeged}, there exists $Q_+ \in \Q, \ N(Q_+)=\St_+$ such that $B_{/ \St_+}=Q_+^*BQ_+.$ Then clearly in this case, the infimum is actually a minimum.
	
	For a non invertible $B,$ consider $\varepsilon >0.$ If $F$ is any lower bound of the set $\{Q_-^*Q_+^*BQ_+Q_- : Q_+ \in \Q, \ N(Q_+)=\St_+\}$ then, for any $Q_+ \in \Q, \ N(Q_+)=\St_+,$ $$F \leq Q_-^*Q_+^*BQ_+Q_- \leq  Q_-^*Q_+^*(B +\varepsilon I)Q_+Q_-.$$ Since $B +\varepsilon I$ is invertible, it follows that $F \leq Q_-^* (B+\varepsilon I)_{/ \St_+} Q_-.$ As $\varepsilon$ is arbitrary, \cite[Corollary 2]{Shorted2} yields $F \leq Q_-^* B_{/ \St_+} Q_-.$
	
\end{dem}

{\vspace{1ex}\noindent{\it Proof of Theorem \ref{Shortedsupinf}.}\hspace{0.5em}} Suppose that $B$ is $\St$-weakly complementable and write $B=B_1+B_2-B_3,$ with $\St \subseteq N(B_1),$ $\St_- \subseteq N(B_2),$ $\St_+ \subseteq N(B_3)$ and $B_2, B_3 \geq_{\HH} 0$ (see Proposition \ref{PropWC3}). Let $Q_- \in \Q,\  N(Q_-)=\St_-.$ Then, for any $Q_+ \in \Q, \ N(Q_+)=\St_+,$  $$Q_-^*Q_+^*BQ_+Q_-= B_1+Q_-^*Q_+^*B_2Q_+Q_- - Q_-^*B_3Q_-.$$ 
By Lemma \ref{lemasupinf}, $\underset{ Q_+ \in \Q, \ N(Q_+)=\St_+}{inf} \  Q_-^*Q_+^*B_2Q_+Q_-=Q_-^*{B_2}_{/ \St_+}Q_-.$ Therefore, 
\begin{align*}
\underset{ Q_+ \in \Q, \ N(Q_+)=\St_+}{inf} \  Q_-^*Q_+^*BQ_+Q_-&=B_1+Q_-^*{B_2}_{/ \St_+}Q_- - Q_-^*B_3Q_-\\
&=B_1+{B_2}_{/ \St_+}- Q_-^*B_3Q_-,
\end{align*}
where we used the fact that $R({B_2}_{/ \St_+}) \subseteq \ol{R(B_2)} \subseteq \St_-^{\perp}$ (see \cite[Corollary 4]{Shorted2}). Finally, by \cite[Proposition 3.7]{Massey} and Lemma \ref{ShortedM1},

$$\underset{ Q_- \in \Q,\  N(Q_-)=\St_-}{sup} \left( \underset{ Q_+ \in \Q, \ N(Q_+)=\St_+}{inf} \  Q_-^*Q_+^*BQ_+Q_-\right) =$$
$$=\underset{ Q_- \in \Q,\  N(Q_-)=\St_-}{sup} \left( B_1+{B_2}_{/ \St_+}- Q_-^*B_3Q_- \right)= B_1+{B_2}_{ / \St_+}-{B_3}_{ / \St_-}=B_{ / \St}.$$

The second equality follows in a similar way.

Conversely, suppose that $$\underset{ Q_- \in \Q,\  N(Q_-)=\St_-}{sup}  \left( \underset{ Q_+ \in \Q, \ N(Q_+)=\St_+}{inf} \  Q_-^*Q_+^*BQ_+Q_- \right)$$ and   $$\underset{ Q_+ \in \Q, \ N(Q_+)=\St_+}{inf} \left( \underset{ Q_- \in \Q,\  N(Q_-)=\St_-}{sup} Q_+^*Q_-^*BQ_-Q_+\right).$$
exist. Then, for every $Q_- \in \Q, N(Q_-)=\St_-,$ there exists $$\underset{Q_+ \in \Q, \ N(Q_+)=\St_+}{inf} \  Q_-^*Q_+^*BQ_+Q_-.$$ In particular, for $Q_-=P_{\St_-^{\perp}},$ $T_0:=\underset{ Q_+ \in \Q, \ N(Q_+)=\St_+}{inf} \ P_{\St_-^{\perp}}Q_+^*BQ_+P_{\St_-^{\perp}}.$ 
Then $P_{\St_-^{\perp}}Q_+^*BQ_+P_{\St_-^{\perp}} - T_0 \geq_{\HH} 0$ for every $Q_+ \in \Q$ such that $N(Q_+)=\St_+.$  
Thus $T_0=P_{\St_-^{\perp}}Q_+^*BQ_+P_{\St_-^{\perp}}-(P_{\St_-^{\perp}}Q_+^*BQ_+P_{\St_-^{\perp}} - T_0 )$ is selfadjoint.

Since $R(Q_+^*BQ) \subseteq \St_+^{\perp},$ then
$P_{\St_-^{\perp}}Q_+^*BQ_+P_{\St_-^{\perp}}=P_{\St^{\perp}}Q_+^*BQ_+P_{\St^{\perp}}$ and $$T_0=\underset{ Q_+ \in \Q, \ N(Q_+)=\St_+}{inf} \ P_{\St^{\perp}}Q_+^*BQ_+P_{\St^{\perp}}.$$ 
Suppose that the matrix decomposition of $T_0$ is given by $T_0= \begin{bmatrix}
t_{11} & t_{12} & t_{13} \\ 
t_{12}^* & t_{22 }& t_{23} \\
t_{13}^* & t_{23}^* &  t_{33} \\ 
\end{bmatrix} \mm{\St_+}{\St_-}{\St^{\perp}}.$ Then 
$$P_{\St^{\perp}}Q_+^*BQ_+P_{\St^{\perp}}-T_0=
\begin{bmatrix}
-t_{11} & -t_{12} & -t_{13} \\ 
-t_{12}^* & -t_{22 }& -t_{23} \\
-t_{13}^* & t_{23}^* &  P_{\St^{\perp}}Q_+^*BQ_+P_{\St^{\perp}}-t_{33} \\ 
\end{bmatrix}  \geq_{\HH} 0.$$ Then $t_{11} \leq_{\HH} 0$ and $t_{22} \leq_{\HH} 0.$ 
Also, since $P_{\St^{\perp}}T_0P_{\St^{\perp}} \leq_{\HH} P_{\St^{\perp}}Q_+^*BQ_+P_{\St^{\perp}}$ for every $Q_+ \in \Q$ such that $N(Q_+)=\St_+,$ then $P_{\St^{\perp}}T_0P_{\St^{\perp}} \leq_{\HH} T_0.$ 
Therefore
$$T_0-P_{\St^{\perp}}T_0P_{\St^{\perp}}=
\begin{bmatrix}
t_{11} & t_{12} & t_{13} \\ 
t_{12}^* & t_{22 }& t_{23} \\
t_{13}^* & t_{23}^* &  0 \\ 
\end{bmatrix}  \geq_{\HH} 0.$$ 
Then, by Lemma \ref{LemmaPositive}, $t_{11} \geq_{\HH} 0$ so $t_{11}=0,$  $t_{22} \geq_{\HH} 0$ so $t_{22}=0,$ and also $t_{12}=t_{12}^*=t_{13}=t_{13}^*=t_{23}=t_{23}^*=0.$ Hence, $R(T_0) \subseteq \St^{\perp}.$ Therefore
$$P_{\St^{\perp}}Q_+^*(B-T_0)Q_+P_{\St^{\perp}} \geq_{\HH}0 \mbox{ for every } Q_+ \in \Q \mbox{ such that } N(Q_+)=\St_+.$$ Let us show that $\PI{(B-T_0)x}{x} \geq 0$ for every $x \in \St_-^{\perp}=\St_+\oplus \St^{\perp}.$
Fix $x \in \St_-^{\perp};$ if $x \in \St_+$ then $\PI{(B-T_0)x}{x}=\PI{Bx}{x} \geq 0,$ because $\St_+ \subseteq N(T_0)$ and $\St_+$ is $B$-nonnegative. If $x \not \in \St_+$ then $P_{\St^{\perp}}x \not = 0$ and there exists a subspace $\M$ such that $x \in \M$ and $\M \dotplus \St_+=\HH.$ Take $Q_+=P_{\M \ {\mathbin{\!/\mkern-3mu/\!}} \ \St_+};$ then $x=Q_+x=Q_+P_{\St^{\perp}}x.$ Thus 
$\PI{(B-T_0)x}{x}=\PI{P_{\St^{\perp}}Q_+^*(B-T_0)Q_+P_{\St^{\perp}}x}{x} \geq0.$
Since $x \in  \St_-^{\perp}$ is arbitrary, $\PI{(B-T_0)x}{x} \geq 0$ for every $x \in \St_-^{\perp}.$ If the matrix decomposition of $B$ is as in \eqref{BS+}, $$P_{\St_-^{\perp}}(B-T_0)P_{\St_-^{\perp}}=\begin{bmatrix}
a_+ & 0 & b_+ \\ 
0& 0 & 0 \\
b_{+}^* & 0 &  c-t_{33} \\ 
\end{bmatrix}  \geq_{\HH} 0.$$ Then, by Lemma \ref{LemmaPositive}, $R(b_+) \subseteq R(a_+^{1/2}).$ 

Analogously, since  $\underset{ Q_+ \in \Q, \ N(Q_+)=\St_+}{inf} \left( \underset{ Q_- \in \Q,\  N(Q_-)=\St_-}{sup} Q_+^*Q_-^*BQ_-Q_+\right)$ exists, it follows that $R(b_-) \subseteq R(a_-^{1/2}).$ Therefore $R(b) \subseteq R(b_+) + R(b_-) \subseteq R(a_+^{1/2}) \oplus R(a_-^{1/2})=R(\vert a \vert^{1/2})$ and $B$ is $\St$-weakly complementable.

{\hfill\qed\vspace{1ex}}

\begin{cor} \label{corpinf} Let $B \in L(\HH)^s$ and let $\St$ be a closed subspace of $\HH.$ Suppose that $\St$ is $B$-nonnegative.
	Then $B$ is $\St$-weakly complementable if and only if there exists $inf \ \{ Q^{*}BQ : Q \in \Q, \ N(Q)=\St \}.$ In this case,
	\begin{equation} \label{infimum}
	B_{/ \St}=inf \ \{ Q^{*}BQ : Q \in \Q, \ N(Q)=\St \}.
	\end{equation}	 
\end{cor}
A similar result holds when $\St$ is $B$-nonpositive, replacing inf by sup.

The following proposition shows that the infimum in \eqref{infimum} is indeed a minimum if and only if $B$ is $\St$-complementable and $\St$ is $B$-nonnegative. 
Recall that when $B$ is $\St$-complementable, the set
$$\mc{P}(B,\St):=\{ Q \in \Q: R(Q)=\St, \ BQ=Q^{*}B \}$$ is not empty. See definition on page 5.

\begin{prop}  \label{Prop36B}   Let $B \in L(\HH)^s$ and let $\St$ be a closed subspace of $\HH.$ Suppose that $B$ is $\St$-weakly complementable. Then
	$$B_{/ \St}=min \ \{ Q^{*}BQ : Q \in \Q, \ N(Q)=\St \}$$ if and only if $B$ is $\St$-complementable and $\St$ is $B$-nonnegative. In this case, 
	$$B_{/ \St}=B(I-Q),$$ for any $Q \in \mc{P}(B,\St).$
\end{prop}

A similar result holds when $\St$ is $B$-nonpositive, replacing min by max.

\begin{dem} Suppose that $B$ is $\St$-complementable and $\St$ is $B$-nonnegative. Then, by \cite[Proposition 4.6]{Massey}, $B_{/ \St}=min \ \{ Q^{*}BQ : Q \in \Q, \ N(Q)=\St \}.$ In this case $B_{/ \St}=B(I-Q),$ for any $Q \in \mc{P}(B,\St).$
	
	Conversely, suppose that $B_{/ \St}=min \ \{ Q^{*}BQ : Q \in \Q, \ N(Q)=\St \}.$ Let $B$ be as in \eqref{matrixW} and $Q:=\matriz{0}{e}{0}{I}$ with $e\in L(\St^{\perp},\St).$ Then $Q \in \Q$ and N$(Q)=\St.$ Let $f$ be the reduced solution of $b=\vert a\vert^{1/2}x$ and  $a=u \vert a \vert$ the polar decomposition of $a.$ From $B_{/ \St} \leq_{\HH} Q^{*}BQ,$ it is easy to check that
	$$0 \leq \PI{u(f+\vert a \vert^{1/2}ue)y)}{(f+\vert a \vert^{1/2}ue)y} \mbox{ for every } y \in \St^{\perp} \mbox{ and } e\in L(\St^{\perp},\St).$$
	Since $R(f)\subseteq \ol{R(\vert a \vert^{1/2})}$ then $\ol{R(\vert a \vert^{1/2})}=\ol{R(f)} \oplus \ol{R(\vert a \vert^{1/2})} \cap R(f)^{\perp}.$ In particular, for every $s \in \St,$ $\vert a \vert^{1/2}s=t+v,$ with $t \in \ol{R(f)}$ and $v \in \ol{R(\vert a \vert^{1/2})} \cap R(f)^{\perp}.$
	If $s\in\St$ and $\vareps >0,$ then there exist $y_{\vareps} \in \St^{\perp}$ and $e_{\vareps} \in L(\St^{\perp},\St)$ such that $\Vert \vert a \vert^{1/2} s - (fy_{\vareps}+\vert a \vert^{1/2}ue_{\vareps}y_{\vareps}) \Vert < \vareps.$ Therefore 
	\begin{align*}
	\PI{as}{s} &=\PI{u\vert a \vert^{1/2} s}{\vert a \vert^{1/2} s}\\
	&= \PI{u \ \underset{\vareps \ra 0}{lim}[(f+\vert a \vert^{1/2}ue_{\vareps})y_{\vareps}]}{\underset{\vareps \ra 0}{lim}[(f+\vert a \vert^{1/2}ue_{\vareps})y_{\vareps}]}\\
	&=\underset{\vareps \ra 0}{lim} \PI{u(f+\vert a \vert^{1/2}ue_{\vareps})y_{\vareps}}{(f+\vert a \vert^{1/2}ue_{\vareps})y_{\vareps}} \geq 0
	\end{align*}
	and $\St$ is $B$-nonnegative. 
	
	In this case, by Proposition \ref{PropWC1}, $B_{/ \St} \in \M^{-}(B,\St^{\perp}),$ so that $B_{/ \St} \leq_{\HH} B.$ Let $Q_0 \in  \Q$ with $N(Q_0)=\St$ such that $B_{/ \St}=Q_0^*BQ_0.$ Then, $B-Q_0^*BQ_0 \geq_{\HH} 0.$ From $Q_0^{*}(B-Q_0^*BQ_0)Q_0=0,$ it follows that $(B-Q_0^*BQ_0)Q_0=0,$ which implies that
	$BQ_0=Q_0^*BQ_0.$ Thus, $E_0:=I-Q_0 \in \mc{P}(B,\St)$ and $B$ is $\St$-complementable.
\end{dem}

If $B \in L(\HH)^s$ and $\St$ is a closed subspace of $\HH$ such that $B$ is $\St$-complementable, then there exists a $B$-selfadjoint projection $Q$ into $\St$ that can be decomposed as the sum of two $B$-selfadjoint projections $Q_+, Q_-$ with $B$-nonnegative and $B$-nonpositive ranges, respectively. 

\begin{lema} \label{THM37} Let $B \in L(\HH)^s$ and let $\St$ be a closed subspace of $\HH.$ Suppose that $\St = \St_{+} \ \oplus_{B} \ \St_{-}$ is any decomposition as in \eqref{WdecompS}. Then, the following statements are equivalent:
	\begin{enumerate}
		\item[i)] $B$ is $\St$-complementable;
		\item [ii)] $B$ is $\St_{\pm}$-complementable;
		\item [iii)] $B$ is $\St$-weakly complementable and $B_{/ \St_{\pm}}$ is $\St_{\mp}$-complementable.
	\end{enumerate}
	In this case, there exists $Q \in \mc{P}(B,\St)$ that can be decomposed as $Q=Q_+ + Q_-,$ where $Q_{\pm} \in \mc{P}(B,\St_{\pm}).$ Moreover, $R(Q_+) \perp R(Q_-)$  and $Q_+ Q_- =Q_-Q_+= 0.$ 
\end{lema}

\begin{dem}  $i) \Leftrightarrow ii):$ Let $\St = \St_{+} \ \oplus_{B} \ \St_{-}$ be any decomposition as in \eqref{WdecompS}. Suppose that the matrix representation of $B$ is as in \eqref{BS+} and $B$ is $\St$-complementable. Then $R(b_+) \subseteq R(a_+)$ and $R(b_-) \subseteq R(a_-).$ In fact, since $R(b) \subseteq R(a),$  for every $y \in \St^{\perp},$ there exists $s \in \St$ such that $by=a s.$ Therefore, for every $y \in \St^{\perp},$ $b_{\pm}y=P_{\St_{\pm}}by = P_{\St_{\pm}} a s=a_{\pm}s$ and $B$ is $\St_{\pm}$-complementable. The converse follows in a similar way using that $R(a)=R(a_+)\oplus R(a_-).$
	
	$i) \Leftrightarrow iii):$ It can be proven in a similar way as in $i) \Leftrightarrow ii)$ using the decomposition of $B_{/ \St_{\pm}}$ given in the proof of Theorem \ref{thmminmax}.
	In this case, let $f$ be the reduced solution of $b_+ = a_+x$ and $g$ the reduced solution of $b_-=- a_-x.$ Set $Q_+:=\begin{bmatrix}
	I & 0 & f  \\ 
	0 & 0 & 0\\
	0  & 0 & 0\\ 
	\end{bmatrix} \mm{\St_+}{\St_-}{\St^{\perp}}$ and $Q_-:=\begin{bmatrix}
	0 & 0 & 0  \\ 
	0 & I & -g\\
	0  & 0 & 0\\ 
	\end{bmatrix} \mm{\St_+}{\St_-}{\St^{\perp}}.$ Then $Q_{\pm} \in \mc{P}(B, \St_{\pm}),$ $R(Q_+) \perp R(Q_-)$ and $Q_+ Q_- =Q_-Q_+= 0.$ Finally, since 
	$$R(f-g) \subseteq \ol{R(a_+)} \oplus \ol{R(a_-)} =\ol{R(a}),$$ $y:= f-g$ is the reduced solution of the equation $b= a x.$ Therefore $Q:=Q_+ + Q_- =\begin{bmatrix}
	I & 0 & f  \\ 
	0 & I & -g\\
	0  & 0 & 0\\ 
	\end{bmatrix} \mm{\St_+}{\St_-}{\St^{\perp}}=\begin{bmatrix}
	I & y \\ 
	0 & 0 \\
	\end{bmatrix} \ms{\St}{\St^{\perp}} \in \mc{P}(B,\St).$
\end{dem}

\begin{cor} \label{Cor310} Let $B \in L(\HH)^s$ and let $\St$ be a closed subspace of $\HH$ such that $B$ is $\St$-complementable.
	Then $$B_{/ \St}=B(I-Q)=(I-Q)^*B$$ for any $Q \in \mc{P}(B,\St).$
\end{cor}

\begin{dem} Let $Q^0 \in \mc{P}(B,\St)$ such that there exist $Q^0_{\pm} \in \mc{P}(B,\St_{\pm})$ with $Q^0=Q^0_+ + Q^0_-,$ $R(Q^0_+) \perp R(Q^0_-)$ and $Q^0_+ Q^0_- =Q^0_-Q^0_+= 0,$ as in Lemma \ref{THM37}. Set $\St_{\pm}:=R(Q^0_{\pm}),$ then $\St=\St_+ \oplus_{B} \St_-,$ as in Lemma \ref{lemaWdecom}.
	By Lemma \ref{THM37}, $B_{ / \St_+}$ is $\St_-$-complementable and since $\St_-$ is $B_{ / \St_+}$-nonpositive, Corollary \ref{corminmax} together with Proposition \ref{Prop36B} give
	$$B_{/ \St}=(B_{/ \St_+})_{ / \St_-}=B_{/ \St_+}(I-Q^0_-).$$ 
	Then, once again by Lemma \ref{THM37}, $B$ is $\St_+$-complementable and by Proposition \ref{Prop36B},
	$$B_{/ \St}=B_{/ \St_+}(I-Q^0_-)=B(I-Q^0_+)(I-Q^0_-)=B(I-Q^0_+-Q^0_-)=B(I-Q^0).$$
	
	Now take any $Q \in \mc{P}(B,\St),$ then by \cite[Theorem 3.5]{CMSSzeged} and \cite[Proposition 3.2]{KreinSzeged}, $Q=Q^0+T,$ for some $T \in L(\HH)$ with $R(T) \subseteq N(B) \cap \St$ and $\St \subseteq N(T).$ Therefore
	$$B_{/ \St}=B(I-Q^0)=B(I-(Q-T))=B(I-Q).$$
\end{dem}

Corollary \ref{Cor310} shows that $B_{/ \St}$ coincides with the Schur complement defined in \cite{GiribetKreinShorted} for a bounded selfadjoint operator $B$  and $\St \subseteq \HH$ a closed subspace such that $B$ is $\St$-complementable.

\begin{cor} \label{corminmaxcomplementable} Let $B \in L(\HH)^s$ and let $\St$ be a closed subspace of $\HH$ such that $B$ is $\St$-weakly complementable. Suppose that $\St = \St_{+} \ \oplus_{B} \ \St_{-}$ is any decomposition as in \eqref{WdecompS}. Then, $B$ is $\St$-complementable if and only if
	\begin{align*}
	B_{/ \St} &=\underset{ Q_- \in \Q,\  N(Q_-)=\St_-}{max}  \left( \underset{ Q_+ \in \Q, \ N(Q_+)=\St_+}{min} \  Q_-^*Q_+^*BQ_+Q_- \right)\\
	&=\underset{ Q_+ \in \Q, \ N(Q_+)=\St_+}{min} \left( \underset{ Q_- \in \Q,\  N(Q_-)=\St_-}{max} Q_+^*Q_-^*BQ_-Q_+\right).
	\end{align*}
\end{cor}

\begin{dem} Suppose that $B$ is $\St$-complementable. Then, by Theorem \ref{Shortedsupinf}, Lemma \ref{THM37} and Corollary \ref{Cor310}, the result follows. 
	
	Conversely, suppose that \begin{align*}
	B_{/ \St} &=\underset{ Q_- \in \Q,\  N(Q_-)=\St_-}{max}  \left( \underset{ Q_+ \in \Q, \ N(Q_+)=\St_+}{min} \  Q_-^*Q_+^*BQ_+Q_- \right)\\
	&=\underset{ Q_+ \in \Q, \ N(Q_+)=\St_+}{min} \left( \underset{ Q_- \in \Q,\  N(Q_-)=\St_-}{max} Q_+^*Q_-^*BQ_-Q_+\right).
	\end{align*} Suppose that the matrix representation of $B$ is as in \eqref{BS+}. Since $B$ is $\St$-weakly complementable, take $B_1, B_2$ and $B_3$ as in the proof of Proposition \ref{PropWC3}. Let $Q_- \in \Q$ with  $N(Q_-)=\St_-.$ Then, by the proof of Theorem \ref{Shortedsupinf},
	$$\underset{ Q_+ \in \Q, \ N(Q_+)=\St_+}{min} \  Q_-^*Q_+^*BQ_+Q_-=Q_-^*(B_1+ {B_2}_{/ \St_+} - B_3)Q_-.$$
	Observe that $$B_1+ {B_2}_{/ \St_+} - B_3= \begin{bmatrix}
	0 & 0 & 0 \\ 
	0 & -a_- & b_- \\
	0 & b_-^*& c-f^*f\\ 
	\end{bmatrix} \mm{\St_+}{\St_-}{\St^{\perp}}.$$ Since there exists $\underset{ Q_- \in \Q,\  N(Q_-)=\St_-}{max} Q_-^*(B_1+ {B_2}_{/ \St_+} - B_3)Q_-,$ by Proposition \ref{Prop36B}, $B_1+ {B_2}_{/ \St_+} - B_3$ is $\St_-$-complementable and then $R(b_-) \subseteq R(a_-).$ In a similar fashion, $B_1+ B_2 - {B_3}_{/ \St_-}$ is $\St_+$-complementable and $R(b_+) \subseteq R(a_+).$ Therefore, $R(b) \subseteq R(b_-) +R(b_+) \subseteq R(a_-) \oplus R(a_+)=R(a).$ Hence $B$ is $\St$-complementable. 
	
\end{dem}

\subsection{A formula for the Schur complement}

When the operator $B$ is $\St$-complementable, the Schur complement can be written as $B_{/ \St}=(I-F)B,$ for any bounded projection with $N(F)=\St^{\perp}$ such that $(FB)^*=FB.$ In fact, from Corollary \ref{Cor310}, it suffices to take $F=Q^*,$ for any  $Q \in \mc{P}(B,\St).$

In this section, we show that a similar formula for $B_{/ \St}$ can be given when $B$ is $\St$-weakly complementable. In this case the projection need not be bounded, but it is densely defined. 

\begin{thm} \label{ThmWC5} Let $B \in L(\HH)^s$ and let $\St$ be a closed subspace of $\HH$. Then $B$ is $\St$-weakly complementable if and only if there exists a densely defined projection $E$ with $N(E)=\St^{\perp}$ such that $E \vert P_{\St}BP_{\St}\vert^ {1/2} \in L(\HH)$ and $EB \in L(\HH)^s.$
	In this case,
	$$B_{ / \St} = (I-E)B.$$ 
\end{thm}

\begin{obs} The densely defined projection $E$ is closed if and only if the pair $(B,\St)$ is quasi-compatible; i.e., $\HH=\ol{\St+ B^{-1}(\St^{\perp})}.$
	Moreover, $E \in L(\HH)$ if and only if  $B$ is $\St$-complementable.
\end{obs}

\begin{dem} Suppose that $B$ is $\St$-weakly complementable. If the matrix decomposition of $B$ induced by $\St$ is as in \eqref{matrixW}, let $f$ be the reduced solution of $b=\vert a\vert^{1/2} x$ and $a=u \vert a \vert$ the polar decomposition of $a.$
	Write $(\vert a \vert^{1/2})^{\dagger}$ for the Moore-Penrose inverse of $\vert a \vert^{1/2}$ and set $$E=\begin{bmatrix} 
	I & 0 \\ 
	f^*u(\vert a \vert^{1/2})^{\dagger} & 0\\  
	\end{bmatrix}.$$
	
	Then $Dom(E)=Dom(\vert a \vert^{1/2})^{\dagger}) \oplus \St^{\perp}$ and $E$ is a densely defined projection with $N(E)=\St^{\perp}.$ On the other hand, since $R(B) \subseteq R(\vert a \vert^ {1/2}) \oplus \St^{\perp},$ the product $(I-E)B$ is well defined. Moreover	
	\begin{align*}
	(I-E)B&=\begin{bmatrix} 
	0 & 0 \\ 
	-f^*u(\vert a \vert^{1/2})^{\dagger} & I\\  
	\end{bmatrix} \begin{bmatrix} 
	\vert a\vert^{1/2} u \vert a\vert^{1/2} & \vert a\vert^{1/2}f \\ 
	f^*\vert a\vert^{1/2} &c\\  
	\end{bmatrix} \\
	& = \begin{bmatrix} 
	0 & 0 \\ 
	0 &c-f^*uf\\  
	\end{bmatrix} =B_{ / \St}
	\end{align*}
	is bounded and selfadjoint. 
	Finally, $$E \vert P_{\St}BP_{\St}\vert^ {1/2}=\begin{bmatrix} 
	\vert a \vert^{1/2} & 0 \\ 
	f^*u & 0\\  
	\end{bmatrix} \in L(\HH).$$

	Conversely, suppose that there exists a densely defined projection $E$ with $N(E)=\St^{\perp}$ such that $E \vert P_{\St}BP_{\St}\vert^ {1/2}$ and $EB \in L(\HH)^s.$ Then the matrix decomposition of $I-E$ is
	$$I-E=\begin{bmatrix} 
	0 & 0 \\ 
	y & I\\  
	\end{bmatrix},$$ with $y: Dom(y) \subseteq \St \ra  \St^{\perp}$ and $\ol{Dom(y)}=\St.$ 
	If the matrix decomposition of $B$ is as in \eqref{matrixW}, since $(I-E)B$ is selfadjoint, it follows that, 
	$ya=-b^*$ and $yb$ is bounded and selfadjoint. 
	From the fact that $E \vert P_{\St}BP_{\St}\vert^ {1/2} \in L(\HH),$ we have that $y\vert a \vert^{1/2} \in L(\St, \St^{\perp})$ and, since $ya=-b^*,$ we also have that $y \vert a \vert^{1/2} u \vert a \vert^{1/2}= -b^*.$ Then $b=\vert a \vert^{1/2} (-y \vert a \vert^{1/2} u )^*,$ $R(b) \subseteq R(\vert a \vert^{1/2})$ and $B$ is $\St$-weakly complementable. 
\end{dem}

\section{Schur complement in Krein spaces}

In this section we adapt the definitions of complementability and weak complementability given in Section \ref{sec3} to a bounded selfadjoint operator $W$ acting on a Krein space $(\HH, \K{ \ } { \ })$ . From now on all spaces are assumed to be Krein spaces unless otherwise stated.

\begin{Def}  Let $W \in L(\HH)^{[s]}$ and let $\St$ be a closed subspace of $\HH.$ The operator $W$ is called $\St$-\emph{complementable} if 
	$$\HH=\St + W^{-1}(\St^{\perpi}).$$
\end{Def}

If $W$ is $\St$-complementable then, for any fundamental decomposition $\HH=\HH_+ \ [\dotplus] \ \HH_-$ with signature operator $J,$ we get that $\HH=\St + (JW)^{-1}(\St^{\perp}).$ Therefore, $W$ is $\St$-complementable if and only if the selfadjoint operator $JW$ is $\St$-complementable in (the Hilbert space) $(\HH, \PI{ \ }{ \ })$ for any (and then for every) signature operator $J.$ From this, it follows that $W$ is $\St$-complementable if and only if there exists a projection $Q$ onto $\St$ such that $WQ=Q^{\#}W.$

In this case, if the matrix representation of $JW$ induced by $\St$ is
\begin{equation} \label{Wdes} 
JW=\begin{bmatrix} 
a & b \\ 
b^* & c \\  
\end{bmatrix}, \end{equation} the $\St$-complementability of $W$ is equivalent to $R(b)\subseteq R(a)$ (see \cite[Proposition 3.3]{CMSSzeged}).

In a similar fashion we define the $\St$-weak complementability in Krein spaces, with respect to a fixed signature operator $J.$

\begin{Def} Let $W \in L(\HH)^{[s]}$ and let $\St$ be a closed subspace of $\HH.$ The operator $W$ is $\St$-\emph{weakly complementable} with respect to a signature operator $J$ if $JW$ is $\St$-weakly complementable in $(\HH, \PI{ \ }{ \ }).$ 
\end{Def}

Next, we show that the $\St$-weak complementability of $W$ does not depend on the signature operator. In order to do so, we need to establish some technical lemmas. Some additional notation is also required: consider $J$ and $J_{\alpha}$ two signature operators and set $\alpha=J_{\alpha}J.$ Denote by $\HH=(\HH, \PI{ \ }{ \ })$ the Hilbert space associated to $J$ and by $\HH_{\alpha}=(\HH, \PI{ \ }{ \ }_{J_\alpha}),$ where $\PI{x }{ y }_{J_\alpha}=\K{J_{\alpha}x}{y}=\PI{\alpha^{-1}x}{y},$ the Hilbert space associated to $J_{\alpha}.$ Then  $\alpha \geq_{\HH} 0$ and $\alpha \geq_{\HH_{\alpha}}0.$  
Notice that $\St^{\perp_{\alpha}}=\alpha (\St^{\perp})$ is the orthogonal complement of $\St$ in $\HH_{\alpha}$ and, for $T \in L(\HH),$ $T^{*_{\alpha}}=\alpha T^* \alpha^{-1}$ is the adjoint of $T$ in $\HH_{\alpha}.$ Also, $T^{*_{\alpha}}=T$ if and only if $\alpha T^*=T \alpha.$ 

Denote by $\vert T \vert_{\alpha}$ the modulus of $T$ in $\HH_{\alpha}$ and by $P_{\St}^{\alpha}=P_{\St {\mathbin{\!/\mkern-3mu/\!}} \St^{\perp_{\alpha}}}$ the orthogonal projection onto $\St$ in $\HH_{\alpha}.$
If $T \geq_{\HH_{\alpha}} 0,$  we indicate by $T^{1/2_{\alpha}}$ the square root of $T$ in $\HH_{\alpha}.$ Frequently, we will use that if $T \geq_{\HH} 0$ or $T=T^*$ then $\alpha T \geq_{\HH_{\alpha}} 0$ or $\alpha T=(\alpha T)^{*_{\alpha}},$ respectively.

\begin{lema} \label{Propalpha0} Let $\St$ be a closed subspace of $\HH,$ $J$ and $J_{\alpha}$ two signature operators and $\alpha=J_{\alpha}J.$ Then
	$$\tilde{\alpha}:=P_{\St}^{\alpha}\alpha |_{\St}=(P_{\St} \alpha^{-1}|_{\St})^{-1} \in GL(\St)^{+}.$$
\end{lema}

\begin{dem} The projection $P_{\St}^{\alpha}$ can be expressed in terms of $P_{\St}$ and $\alpha$ as
	$$P_{\St}^{\alpha}=P_{\St}(P_{\St}\alpha^{-1}P_{\St}+(I-P_{\St}) \alpha^{-1} (I-P_{\St}))^{-1} \alpha^{-1};$$
	see  \cite[Section 4]{CMSSzeged}. Therefore, 
	$$P_{\St}^{\alpha}\alpha |_{\St}=P_{\St}(P_{\St}\alpha^{-1}P_{\St}+(I-P_{\St}) \alpha^{-1} (I-P_{\St}))^{-1} |_{\St}=$$$$=(P_{\St} \alpha^{-1}|_{\St})^{-1}  \in GL(\St)^{+}.$$
\end{dem}

By Lemma \ref{Propalpha0}, $\St_{\tilde{\alpha}}:=(\St, \PI{\cdot}{\cdot}_{\tilde{\alpha}})$ is a Hilbert space. Also, from the discussion before Lemma \ref{Propalpha0}, if $a \in L(\St)^s$ then $(\tilde{\alpha}a)^{*_{\tilde{\alpha}}}=\tilde{\alpha}a.$
\begin{lema} \label{Coralpha1} Let $\St$ be a closed subspace of $\HH,$ $J$ and $J_{\alpha}$ two signature operators, $\alpha=J_{\alpha}J$ and $a\in L(\St)^s.$ Let $a'= P_{\St}^{\alpha} \alpha a.$ Then
	\begin{equation} \label{rangeeq}
	R(\vert a' \vert_{\tilde{\alpha}}^{1/2_{\tilde{\alpha}}})=R((P_{\St}^{\alpha} \alpha \vert a \vert P_{\St}^{\alpha})^{1/2_{\alpha}})=\tilde{\alpha}R(\vert a \vert^{1/2}).
	\end{equation}
	
\end{lema}

\begin{dem} First observe that  $$(P_{\St}^{\alpha} \alpha \vert a \vert P_{\St}^{\alpha})^{1/2_{\alpha}}=(\tilde{\alpha} \vert a \vert P_{\St}^{\alpha})^{1/2_{\alpha}}=P_{\St}^{\alpha} \alpha X_0 P_{\St}^{\alpha},$$ with
	$X_0=\tilde{\alpha}^{-1/2}(\tilde{\alpha}^{1/2}\vert a \vert \tilde{\alpha}^{1/2})^{1/2}\tilde{\alpha}^{-1/2}.$ In fact, $X_0 \in L(\St)^{+},$ because $\tilde{\alpha}^{1/2} \in GL(\St)^{+}$ by Lemma \ref{Propalpha0}. Clearly, 
	$X_0 \tilde{\alpha} X_0=\vert a \vert.$ Therefore, $(\tilde{\alpha}X_0)^2=\tilde{\alpha}\vert a \vert.$ 
	Also, 
	$\tilde{\alpha}X_0 P_{\St}^{\alpha} = P_{\St}^{\alpha} \alpha X_0 P_{\St}^{\alpha} \geq_{\HH_{\alpha}} 0.$ 
	
	Then $(\tilde{\alpha}X_0P_{\St}^{\alpha})^2=\tilde{\alpha}X_0P_{\St}^{\alpha}\tilde{\alpha}X_0P_{\St}^{\alpha}=\tilde{\alpha}X_0\tilde{\alpha}X_0P_{\St}^{\alpha}=\tilde{\alpha}\vert a \vert P_{\St}^{\alpha}.$
	Thus, $$(\tilde{\alpha} \vert a \vert P_{\St}^{\alpha})^{1/2_{\alpha}}=\tilde{\alpha}X_0P_{\St}^{\alpha}=P_{\St}^{\alpha} \alpha X_0 P_{\St}.$$
	
	Now, since $X_0 \tilde{\alpha} X_0=\vert a \vert,$ Douglas' Lemma yields $R(X_0\tilde{\alpha}^{1/2})=R(\vert a \vert^{1/2}).$ Therefore, because $\tilde{\alpha}^{1/2} \in GL(\St)^+$, $R(X_0)=R(\vert a \vert^{1/2})$ (see Lemma \ref{Propalpha0}). Then 
	$$R((P_{\St}^{\alpha} \alpha \vert a \vert P_{\St}^{\alpha})^{1/2_{\alpha}})=R(P_{\St}^{\alpha} \alpha X_0 P_{\St})=R(\tilde{\alpha}X_0P_{\St}^{\alpha})=\tilde{\alpha}R(X_0)=\tilde{\alpha}R(\vert a \vert^{1/2})$$ and the second equality in \eqref{rangeeq} follows.
	Using that $(\tilde{\alpha}a)^{*_{\tilde{\alpha}}}=\tilde{\alpha}a,$
	$R(\vert a' \vert_{\tilde{\alpha}})=R(a')=R(\tilde{\alpha}a)=R(\tilde{\alpha}\vert a \vert).$ Then, applying Douglas' Lemma and the operator monotonicity of the square root in $\St_{\tilde{\alpha}}$ (see \cite{Pedersen}), we get that
	$$R(\vert a' \vert_{\tilde{\alpha}}^{1/2_{\tilde{\alpha}}})=R((\tilde{\alpha} \vert a \vert)^{1/2_{\tilde{\alpha}}}).$$
	Finally, from $R((\tilde{\alpha} \vert a \vert)^{1/2_{\tilde{\alpha}}})=R((P_{\St}^{\alpha} \alpha \vert a \vert P_{\St}^{\alpha})^{1/2_{\alpha}}),$ we get the first equality.
\end{dem}

\begin{lema} \label{Coralpha3} Let $\St$ be a closed subspace of $\HH,$ $J$ and $J_{\alpha}$ two signature operators, $\alpha=J_{\alpha}J$ and $a\in L(\St)^s.$ Let $a'= P_{\St}^{\alpha} \alpha a,$ $\Gamma:=\begin{bmatrix} 
	\vert a \vert^{1/2} & 0 \\ 
	0 & I\\  \end{bmatrix} \ms{\St}{\St^{\perp}},$  $\Gamma_{\alpha}:=\begin{bmatrix} 
	\vert a'  \vert_{\tilde{\alpha}}^{1/2_{\tilde{\alpha}}} & 0 \\ 
	0 & I\\  
	\end{bmatrix} \ms{\St}{\St^{\perp\alpha}}$ and $E$ a densely defined projection with $N(E)=\St^{\perp}$ such that $ E\Gamma \in L(\HH).$ Then $\alpha E \alpha^{-1}$ is a densely defined projection with $N(\alpha E \alpha^{-1})=\St^{\perp_{\alpha}}$ such that $\alpha E \alpha^{-1} \Gamma_{\alpha} \in L(\HH).$ 
\end{lema}

\begin{dem} Clearly, $\alpha E \alpha^{-1}$ is a densely defined projection with $N(\alpha E \alpha^{-1})=\St^{\perp_{\alpha}}.$
	Let us see that $ \alpha E \alpha^{-1} \Gamma_{\alpha} \in L(\HH).$  Since $E$ is a densely defined projection with $N(E)=\St^{\perp}$ such that $E\Gamma \in L(\HH),$  the matrix decomposition of $E$ is
	$$E=\begin{bmatrix} 
	I & 0 \\ 
	y & 0\\  
	\end{bmatrix} \ms{\St}{\St^{\perp}}=(I+y)P_{\St},$$ with $y: Dom(y) \subseteq \St \ra  \St^{\perp},$ $\ol{Dom(y)}=\St$ and $y \vert a \vert^{1/2} \in L(\St, \St^{\perp}).$
	Also, $\Gamma_{\alpha}=\vert a'  \vert_{\tilde{\alpha}}^{1/2_{\tilde{\alpha}}} P_{\St}^{\alpha} + (I-P_{\St}^{\alpha}).$ 
	Since, by Lemma \ref{Coralpha1}, $R(\vert a' \vert_{\tilde{\alpha}}^{1/2_{\tilde{\alpha}}})=R(\tilde{\alpha}\vert a \vert^{1/2}),$ $\vert a'\vert_{\tilde{\alpha}}^{1/2_{\tilde{\alpha}}} \geq_{\St_{\tilde{\alpha}}} 0$ and $\tilde{\alpha}\vert a \vert^{1/2} \geq_{\St_{\tilde{\alpha}}} 0,$ there exists $g \in GL(\St)$ such that $$\vert a' \vert_{\tilde{\alpha}}^{1/2_{\tilde{\alpha}}}=\tilde{\alpha}\vert a \vert^{1/2}g.$$
	Then $ P_{\St} \alpha^{-1} \vert a' \vert_{\tilde{\alpha}}^{1/2_{\tilde{\alpha}}}=P_{\St} \alpha^{-1} \tilde{\alpha}\vert a \vert^{1/2}g= P_{\St} \alpha^{-1}  P_{\St} \tilde{\alpha}\vert a \vert^{1/2}g= \tilde{\alpha}^{-1} \tilde{\alpha} \vert a \vert^{1/2}g=\vert a \vert^{1/2}g.$ 
	
	Therefore, $R(P_{\St} \alpha^{-1} \vert a' \vert_{\tilde{\alpha}}^{1/2_{\tilde{\alpha}}}) = R(\vert a \vert^{1/2}) \subseteq Dom(y)$ and $$yP_{\St} \alpha^{-1} \vert a' \vert_{\tilde{\alpha}}^{1/2_{\tilde{\alpha}}}=y \vert a \vert^{1/2}g \in L(\St, \St^{\perp}).$$
	
	Thus
	$$\alpha E \alpha^{-1} \Gamma_{\alpha}= \alpha P_{\St}\alpha^{-1} \vert a'  \vert_{\tilde{\alpha}}^{1/2_{\tilde{\alpha}}} P_{\St}^{\alpha}+\alpha yP_{\St} \alpha^{-1}\vert a'  \vert_{\tilde{\alpha}}^{1/2_{\tilde{\alpha}}} P_{\St}^{\alpha}  \in L(\HH).$$ 
\end{dem}

Now we are ready to show that the weak $\St$-complementability of $W$ does not depend on the fundamental decomposition of $\HH.$ 

\begin{thm} \label{WC} Let $W \in L(\HH)^{[s]}$ and let $\St$ be a closed subspace of $\HH.$ Suppose that $W$ is $\St$-weakly complementable for some signature operator $J.$ Then $W$ is $\St$-weakly complementable for any other signature operator $J_{\alpha}.$
\end{thm}

\begin{dem} Suppose that $W$ is $\St$-weakly complementable for some signature operator $J$ and the matrix decomposition of $JW$ is as in \eqref{Wdes}. Observe that $a \in L(\St)^s.$ 
	
	By Theorem \ref{ThmWC5},  there exists a densely defined projection $E$ with $N(E)=\St^{\perp}$ such that $EJW \in L(\HH)^s.$ Also, if $\Gamma:=\begin{bmatrix} 
	\vert a \vert^{1/2} & 0 \\ 
	0 & I\\  
	\end{bmatrix} \ms{\St}{\St^{\perp}},$ then $E\Gamma \in L(\HH).$ Let $J_{\alpha}$ be another signature operator and $\alpha=J_{\alpha}J.$ If $J_{\alpha}W=\alpha JW=\begin{bmatrix} 
	a' & b' \\ 
	b'^* & c' \\  
	\end{bmatrix} \ms{\St}{\St^{\perp_{\alpha}}}$ then $$a'= P_{\St}^{\alpha} \alpha JW P_{\St}^{\alpha}=P_{\St}^{\alpha} \alpha a P_{\St}^{\alpha}.$$ Consider $\alpha E \alpha^{-1}$ and $\Gamma_{\alpha}:=\begin{bmatrix} 
	\vert a'  \vert_{\tilde{\alpha}}^{1/2_{\tilde{\alpha}}} & 0 \\ 
	0 & I\\  
	\end{bmatrix} \ms{\St}{\St^{\perp\alpha}}.$ By Lemma \ref{Coralpha3},  $\alpha E \alpha^{-1}$ is a densely defined projection with $N(\alpha E \alpha^{-1})=\St^{\perp_{\alpha}}$ such that $\alpha E \alpha^{-1} \Gamma_{\alpha} \in L(\HH).$ Also, $(\alpha E \alpha^{-1})(J_{\alpha}W)=\alpha EJW$ is selfadjoint in $\HH_{\alpha}.$ Therefore, again by Theorem \ref{ThmWC5}, $J_{\alpha}W$ is $\St$-weakly complementable.
\end{dem}

From now on, since the $\St$-weak complementability does not depend on the fundamental decomposition of $\HH,$ we simply say that $W$ is $\St$-weakly complementable, whenever $W$ is $\St$-weakly complementable with respect to a signature operator $J.$ In particular, if $W \geq0$ then $W$ is $\St$-weakly complementable.

Following the ideas of \cite{AntCorSto06}, we extend the notion of Schur complement to selfadjoint operators in Krein spaces:

\begin{Def} Let $W \in L(\HH)^{[s]},$ $\St$ a closed subspace of $\HH$ and $J$ a signature operator. Suppose that $W$ is $\St$-weakly complementable. The \emph{Schur complement} of $W$ to $\St$ corresponding to $J$ is
	$$W_{/ [\St]}^J =J (JW)_{ / \St},$$
	and the $\St$-\emph{compression} of $W$ is $W_{ [\St]}^J = W- W_{/ [\St]}^J=J (JW)_{\St}.$
\end{Def}

\begin{thm}\label{thethm} Let $W \in L(\HH)^{[s]}$ and let $\St$ be a closed subspace of $\HH.$ Suppose that $W$ is $\St$-weakly complementable and $J$ is a signature operator, then 
	$$W_{/ [\St]}^J =W_{/ [\St]}^{J_{\alpha}},$$
	for any other signature operator $J_{\alpha};$ i.e., the Schur complement does not depend on the fundamental decomposition of $\HH.$ 
\end{thm}
Henceforth we write $W_{/ [\St]}$ for this operator.

\begin{dem} Suppose that $W$ is $\St$-weakly complementable. Then, by Lemma \ref{Coralpha3}, $JW$ is $\St$-weakly complementable for any signature operator $J.$ By Theorem \ref{ThmWC5},  there exists a densely defined projection $E$ with $N(E)=\St^{\perp}$ such that $EJW \in L(\HH)^s,$ $E\Gamma \in L(\HH)$ and  $W_{/ [\St]}^J=J(I-E)JW.$ Let $J_{\alpha}$ be another signature operator, $\alpha=J_{\alpha}J$ and consider $\alpha E \alpha^{-1}.$ Then, by Theorem \ref{WC}, $\alpha E \alpha^{-1}$ is a densely defined projection with $N(\alpha E \alpha^{-1})=\St^{\perp_{\alpha}}$ such that $\alpha E \alpha^{-1} \Gamma_{\alpha} \in L(\HH).$ Therefore, by Theorem \ref{ThmWC5}, $$W_{/ [\St]}^{J_{\alpha}}=J_{\alpha} \alpha( I- E) \alpha^{-1}J_{\alpha}W=J(I-E)JW=W_{/ [\St]}^J.$$ 
\end{dem}

\begin{cor} Let $W \in L(\HH)^{[s]}$ and let $\St$ be a closed subspace of $\HH.$ Suppose that $W$ is $\St$-weakly complementable. Then there exists a densely defined projection $E$ with $N(E)=\St^{\perpi}$ such that
	$$W_{/ [\St]}=(I-E)W.$$
\end{cor}

\bigskip
Let $W \in L(\HH)^{[s]}$ and let $\St$ be a closed subspace of $\HH,$ define
$$\N^{-}(W,\St^{\perpi}):=\{ X \in L(\HH)^{[s]} : X \leq W, \ R(X) \subseteq \St^{\perpi} \},$$
$$\N^{+}(W,\St^{\perpi}):=\{ X \in L(\HH)^{[s]} : W \leq X, \ R(X) \subseteq \St^{\perpi}  \}.$$
If $J$ is any signature operator,
$$\N^{\pm}(W,\St^{\perpi})=J\M^{\pm}(JW,\St^{{\perp}}).$$

Let $W \in L(\HH)^{[s]}$ and let $\St$ be a closed subspace of $\HH.$ Then, applying Lemma \ref{lemaWdecom} to $B=JW,$ with $J$ any signature operator, $\St$ can be decomposed as 
\begin{equation} \label{WdecompSKrein}
\St=\St_+ \ [\dotplus]_{W} \ \St_-,
\end{equation}
where $\St_{+}$ and $\St_-$ are closed, $\St_+$ is $W$-nonnegative, $\St_-$ is $W$-nonpositive and, moreover, $\St_{+} \perp \St_{-}.$

\begin{prop} \label{PropWCK3} Let $W \in L(\HH)^{[s]}$ and let $\St$ be a closed subspace of $\HH.$ Suppose that $\!$ $\St=\St_+ \ [\dotplus]_{W} \ \St_-$ is any decomposition as in \eqref{WdecompSKrein} for some signature operator $J.$ Then the following statements are equivalent:
	\begin{enumerate}
		\item[i)] $W$ is $\St$-weakly complementable;
		\item [ii)] there exist $W_1, W_2, W_3 \in L(\HH)^{[s]},$ $W_2, W_3 \geq 0$ such that $W=W_1+W_2-W_3$ and $\St \subseteq N(W_1),$ $\St_- \subseteq N(W_2),$ $\St_+ \subseteq N(W_3);$  
		\item [iii)] the sets $\N^{-}(W,\St_{+}^{\perpi})$ and $\N^{+}(W,\St_{-}^{\perpi})$ are non-empty;
		\item [iv)] $W$ is $\St_{\pm}$-weakly complementable.
	\end{enumerate}
\end{prop}

\begin{dem} This follows from Proposition \ref{PropWC3}.
\end{dem} 

The following theorem proves that the set $\N^{-}(W,\St^{\perpi})$ has a maximum element if and only if $\St$ is $W$-nonnegative and $W$ is $\St$-weakly complementable. A similar result can be proven if $\St$ is $W$-nonpositive and $W$ is $\St$-weakly complementable.

\begin{prop} \label{ShortedC1} Let $W \in L(\HH)^{[s]}$ and let $\St$ be a closed subspace of $\HH.$  Then $\St$ is $W$-nonnegative and $W$ is $\St$-weakly complementable if and only if the set $\N^{-}(W,\St^{\perpi})$ has a maximum element.
	
	In this case, $$\WS= max \ \N^{-}(W,\St^{\perpi}).$$ 
\end{prop}

\begin{dem}  This follows from Proposition \ref{PropWC1}.
\end{dem}

\begin{thm} \label{thmminmaxII} Let $W \in L(\HH)^{[s]}$ and let $\St$ be a closed subspace of $\HH.$ Suppose that $\St=\St_+ \ [\dotplus]_{W} \ \St_-$ is any decomposition as in \eqref{WdecompSKrein} for some signature operator $J.$ Then $W$ is $\St$-weakly complementable if and only if there exists $min \ \N^{+} (max  \ \N^{-}(W,\St_{+}^{\perpi}), \St_{-}^{\perpi})$ and $max \ \N^{-} (min  \ \N^{+}(B,\St_{-}^{\perpi}), \St_{+}^{\perpi}).$
	
	In this case,
	\begin{align*}
	\WS&= min \ \N^{+} (max  \ \N^{-}(W,\St_{+}^{\perpi}), \St_{-}^{\perpi})\\
	&=max \ \N^{-} (min  \ \N^{+}(B,\St_{-}^{\perpi}), \St_{+}^{\perpi}).
	\end{align*}
	
\end{thm}

\begin{dem} This follows from Theorem \ref{thmminmax}.
\end{dem}

\begin{cor} \label{corminmaxII} Let $W \in L(\HH)^{[s]}$ and let $\St$ be a closed subspace of $\HH.$ Suppose that $\St=\St_+ \ [\dotplus]_{W} \ \St_-$ is any decomposition as in \eqref{WdecompSKrein} for some signature operator $J.$  If $W$ is $\St$-weakly complementable, then 
	$$\WS= (W_{/ [ \St_+]} )_{/ [\St_-]}=(W_{/[ \St_-]})_{/ [\St_+]}.$$
\end{cor}

\begin{thm} \label{ShortedsupinfKrein} Let $W \in L(\HH)^{[s]}$ and let $\St$ be a closed subspace of $\HH.$ Suppose that $\St=\St_+ \ [\dotplus]_{W} \ \St_-$ is any decomposition as in \eqref{WdecompSKrein} for some signature operator $J.$ $W$ is $\St$-weakly complementable if and only if there exist $$\underset{ E_- \in \Q,\  N(E_-)=\St_-}{sup} \left( \underset{ E_+ \in \Q, \ N(E_+)=\St_+}{inf} \  E_-^{\#}E_+^{\#}WE_+E_- \right)$$ and
	$$\underset{ E_+ \in \Q, \ N(E_+)=\St_+}{inf} \left(\underset{ E_- \in \Q,\  N(E_-)=\St_-}{sup} E_+^{\#}E_-^{\#}WE_-E_+\right).$$
	In this case, 
	\begin{align*}
	\WS&=\underset{ E_- \in \Q,\  N(E_-)=\St_-}{sup} \left( \underset{ E_+ \in \Q, \ N(E_+)=\St_+}{inf} \  E_-^{\#}E_+^{\#}WE_+E_- \right) \\
	&=\underset{ E_+ \in \Q, \ N(E_+)=\St_+}{inf} \left(\underset{ E_- \in \Q,\  N(E_-)=\St_-}{sup} E_+^{\#}E_-^{\#}WE_-E_+\right).
	\end{align*}
	
\end{thm}

\begin{dem} For any signature operator $J,$ if $(\HH, \PI{ \ }{ \ })$ is the associated Hilbert space, 
	\begin{equation} \label{JMN}
	J \{ Q^{*}JWQ : Q \in \Q, \ N(Q)=\St \}=\{ E^{\#}WE : E \in \Q, \ N(E)=\St \}.
	\end{equation}
	Also, there exists $inf_{\leq} \ \{ E^{\#}WE: E=E^2, \ N(E)=\St\}$ if and only if there exists $inf_{\leq_{\HH}}  \{ Q^{*}JWQ : Q \in \Q, \ N(Q)=\St\}.$ Moreover 
	\begin{equation} \label{JMN2}
	inf_{\leq} \ \{ E^{\#}WE: E=E^2, \ N(E)=\St\}=Jinf_{\leq_{\HH}} \{ Q^{*}JWQ : Q \in \Q, \ N(Q)=\St\}.
	\end{equation}
	
	Analogously,  there exists $sup_{\leq} \ \{ E^{\#}WE: E=E^2, \ N(E)=\St\}$ if and only if there exists $\linebreak$ $sup_{\leq_{\HH}}  \{ Q^{*}JWQ : Q \in \Q, \ N(Q)=\St\}.$ Moreover 
	\begin{equation} \label{JMN3}
	sup_{\leq} \ \{ E^{\#}WE: E=E^2, \ N(E)=\St\}=Jsup_{\leq_{\HH}} \{ Q^{*}JWQ : Q \in \Q, \ N(Q)=\St\}.
	\end{equation}
	
	The result follows from \eqref{JMN}, \eqref{JMN2}, \eqref{JMN3} and Theorem \ref{Shortedsupinf}.
\end{dem}

\begin{cor} \label{ShortedC2} Let $W \in L(\HH)^{[s]}$ and let $\St$ be a closed subspace of $\HH.$  Suppose that $\St$ is $W$-nonnegative. Then $W$ is $\St$-weakly complementable if and only if there exists $inf \ \{ E^{\#}WE: E=E^2, \ N(E)=\St\}.$
	
	In this case, $$\WS=inf \ \{ E^{\#}WE: E=E^2, \ N(E)=\St\}.$$ 
\end{cor}

A similar result holds when $\St$ is $W$-nonpositive, replacing inf by sup.

\begin{prop}  \label{Prop36} Let $W \in L(\HH)^{[s]}$ and let $\St$ be a closed subspace of $\HH.$ Suppose that $W$ is $\St$-weakly complementable.  Then,
	$$\WS =min \ \{ E^{\#}WE : E \in \Q, \ N(E)=\St \}$$ if and only if $W$ is $\St$-complementable  and $\St$ is $W$-nonnegative.
	
	In this case, 
	$$\WS=W(I-Q),$$ with $Q$ any projection onto $\St$ such that $WQ=Q^{\#}W.$
\end{prop}
A similar result holds when $\St$ is $W$-nonpositive, replacing min by max.

\begin{dem} For any signature operator $J,$ if $(\HH, \PI{ \ }{ \ })$ is the associated Hilbert space, by \eqref{JMN} and Proposition \ref{Prop36B},
	\begin{align*}
	\WS&=J (JW)_{ / \St}=J \ min_{ \leq_{\HH}} \{ Q^{*}JWQ : Q \in \Q, \ N(Q)=\St \}\\
	&= min \ \{ E^{\#}WE : E \in \Q, \ N(E)=\St \}
	\end{align*} if and only if $\HH=\St+ (JW)^{-1}(\St^{\perp})$ and $\St$ is $JW$-nonnegative (in the Hilbert space $\HH$) if and only if $W$ is $\St$-complementable  and $\St$ is $W$-nonnegative.
	
	The operator $Q$ is any projection onto $\St$ such that $WQ=Q^{\#}W$ if and only if $Q \in \mc{P}(JW,\St)$ for any signature operator $J.$ Therefore, in these cases, by Proposition \ref{Prop36B},  
	$\WS=J (JW)_{ / \St}=J (JW)(I-Q)=W(I-Q),$ for any of these projections.
\end{dem}

\begin{cor} \label{Cor310W} Let $W \in L(\HH)^{[s]}$ and let $\St$ be a closed subspace of $\HH$ such that $W$ is $\St$-complementable.
	Then $$\WS=W(I-Q),$$ for $Q$ any projection onto $\St$ such that $WQ=Q^{\#}W.$
\end{cor}

\begin{dem} This follows proceeding as in Corollary \ref{Cor310} and by Proposition \ref{Prop36}.
\end{dem}

\begin{thm} Let $W \in L(\HH)^{[s]}$ and let $\St$ be a closed subspace of $\HH$ such that $W$ is $\St$-weakly complementable. Suppose that $\St=\St_+ \ [\dotplus]_{W} \ \St_-$ is any decomposition as in \eqref{WdecompSKrein} for some signature operator $J$. Then, $W$ is $\St$-complementable if and only if
	\begin{align*}
	W_{/ [\St]}&=\underset{ E_+ \in \Q, \ N(E_+)=\St_+}{min} \left( \underset{ E_- \in \Q,\  N(E_-)=\St_-}{max} E_-^{\#}E_+^{\#}WE_+E_- \right)\\
	&=\underset{ E_- \in \Q,\  N(E_-)=\St_-}{max} \left( \underset{ E_+ \in \Q, \ N(E_+)=\St_+}{min} \ E_+^{\#}E_-^{\#}WE_-E_+  \right).
	\end{align*}
	
\end{thm}
\begin{dem} This follows by \eqref{JMN} and Corollary \ref{corminmaxcomplementable}.
\end{dem}

\subsection{Weak complementability for regular subspaces}

Any $W \in L(\HH)^{[s]}$ can be written in the form $$W=DD^{\#}$$ where $D \in L(\KK,\HH)$ for some Krein space $\KK$ and $N(D)=\{0\}.$  This factorization, in general, is not unique. Such factorizations are known as \emph{Bogn\'ar-Kr\'amli factorizations}, see \cite{BognarKramli}.

Let $J$ be any signature operator of $\HH.$ Then, $JW$ is selfadjoint in the corresponding Hilbert space. If $JW=U \vert JW \vert = \vert JW \vert U$ is the polar factorization of $JW,$ then $\KK:= \ol{R(\vert JW \vert)}$ is a Krein space with signature operator $J_{\KK}:=U|_{\KK}.$ Define $D : \KK \ra \HH$ by 
\begin{equation} \label{Deq}
Dk:=J \vert JW \vert^{1/2}k, \  k \in \KK.
\end{equation}
Then, $N(D)=\{0\},  $ $D^{\#}=J_{\KK} \vert JW \vert^{1/2}=U\vert JW \vert^{1/2}$ and  $DD^{\#}=W$ (cf. \cite[Theorem 1.1]{DR1}).

\begin{Def} A Bogn\'ar-Kr\'amli factorization of an operator $W \in L(\HH)^{[s]}$ which is constructed by the method described above is called a \emph{polar factorization} of $W$ (see \cite[Lecture 6]{DR1}).
\end{Def}

\begin{lema} \label{polarfact} Let $W \in L(\HH)^{[s]}$ have polar factorizations $W=DD^{\#}=EE^{\#}$ with $D\in L(\KK, \HH)$ and $E \in L(\KK',\HH).$ Then
	$$R(D)=R(E).$$
\end{lema}

\begin{dem} In this case, following similar arguments as in \cite[Theorem 6.1]{DR1} and \cite[Theorem 6.2]{DR1}, it can be shown that there exists a unique $L \in L(\KK', \KK)$ such that $E=DL$ and $D=EL^{\#}.$ Clearly, $R(D)=R(E).$
\end{dem}

Let $\St$ be a regular subspace of $\HH,$ then $W \in L(\HH)^{[s]}$ can be represented as a $2\times 2$ block matrix in the form
\begin{equation} \label{WmatrixI}
W=\matriz{w_{11}}{w_{12}}{w_{12}^{\#}}{w_{22}} \ms{\St}{\St^{\perpi}}.
\end{equation}

\begin{thm} \label{ThmShortedKrein} Let $W \in L(\HH)^{[s]}$ and let $\St$ be a regular subspace of $\HH.$ Suppose that $W$ is represented as in \eqref{WmatrixI} and  $w_{11}=dd^{\#}$ is any polar factorization of $w_{11}.$ Then $W$ is $\St$-weakly complementable if and only if $R(w_{12}) \subseteq R(d).$
	
	In this case, $$W_{/ [\St]}=\matriz{0}{0}{0}{w_{22}-y^{\#}y} \ms{\St}{\St^{\perpi}},$$ with $y \in L(\St^{\perpi}, \KK)$ the only solution of the equation $w_{12}=dx.$
\end{thm}

\begin{dem} Take $J=\matriz{J_1}{0}{0}{J_2}  \ms{\St}{\St^{\perpi}}$ a signature operator for $\HH$ as in \eqref{signature}. 
	
	Then, $JW=\matriz{J_1w_{11}}{J_1w_{12}}{J_2w_{12}^{\#}}{J_2w_{22}} \ms{\St}{\St^{\perp}}$ and $W$ is $\St$-weakly complementable if and only if $R(J_1w_{12}) \subseteq R(\vert J_1 w_{11}\vert^ {1/2})$ or, equivalently, $R(w_{12}) \subseteq R(J_1\vert J_1 w_{11}\vert^ {1/2})=R(d).$ Indeed, if  $e:=J_1\vert J_1 w_{11}\vert^ {1/2}$ then, by \eqref{Deq}, $w_{11}=ee^{\#}$ is a polar factorization of $w_{11}$ and, by Lemma \ref{polarfact}, $R(e)=R(d).$ 
	
	In this case, let $y \in L(\St^{\perpi}, \KK)$ be the only solution of the equation $w_{12}=dx.$ Observe that
	$$(JW)_{/ \St}=\matriz{0}{0}{0}{J_2w_{22}-f^*uf} \ms{\St}{\St^{\perp}},$$ where $f$ is the reduced solution of $J_1w_{12}=\vert J_1w_{11} \vert^{1/2} x$ and $u$ is the partial isometry corresponding to the polar decomposition of $J_1w_{11}.$
	Clearly, $w_{12}=J_1 \vert J_1w_{11} \vert^{1/2} f=ef,$ so that $ef=dy.$ As in the  proof of Lemma \ref{polarfact}, a unique bounded operator $l$ can be found such that $e=dl$ and $d=el^{\#}.$ Therefore, $y^{\#}d^{\#}=f^{\#}e^{\#}=f^{\#}l^{\#}d^{\#}$ and, since $d^{\#}$ has a dense range, $y^{\#}=f^{\#}l^{\#}.$ In a similar way, $l^{\#}y=f.$ Therefore, $y^{\#}y=f^{\#}l^{\#}y=f^{\#}f.$ Finally, since $J_2f^*u=f^{\#},$ it follows that
	$$W_{/ [\St]}=J(JW)_{/ \St}=\matriz{0}{0}{0}{w_{22}-J_2f^*uf} =\matriz{0}{0}{0}{w_{22}-y^{\#}y} \ms{\St}{\St^{\perpi}}.$$
\end{dem}

\subsection{An application to a completion problem}

Let $\St$ be a regular subspace of $\HH$ and consider a bounded incomplete block operator
\begin{equation} \label{WmatrixIncomp}
W^0=\matriz{w_{11}}{w_{12}}{w_{12}^{\#}}{*} \ms{\St}{\St^{\perpi}},
\end{equation}
with $w_{11} \in L(\St)^{[s]}.$

Following the ideas of Baidiuk in \cite[Theorem 2.1]{Baidiuk}, the next theorem solves a completion problem for any bounded incomplete operator $W^0$ of the form \eqref{WmatrixIncomp}.

\begin{prop} \label{PropComple} Let $\St$ be a regular subspace of $\HH$ and $W^0$ be an incomplete block operator of the form \eqref{WmatrixIncomp}. Assume that the number of negative squares $\nu_-[w_{11}]$ of the quadratic form $\K{w_{11}f}{f}, \ f \in \St,$ is finite. Let $w_{11}=dd^{\#}$ be any polar factorization of $w_{11}.$ 
	Then, there exists a completion $W$ of $W^0$ with some operator $w_{22} \in L(\St^{\perpi})^{[s]}$ such that $\nu_-[W]=\nu_-[w_{11}]$ if and only if $R(w_{12}) \subseteq R(d).$
	
	In this case, if $y$ is the unique bounded solution of the equation $w_{12}=dx,$ the operator $y^{\#}y \in L(\St^{\perpi})$ is the minimum in the solution set
	$$\mc{W}= \{ w_{22} \in L(\St^ {\perpi})^{[s]}: W=\matriz{w_{11}}{w_{12}}{w_{12}^{\#}}{w_{22}} : \nu_-[W]=\nu_-[w_{11}] \},$$ and this solution set admits the description 
	$$\mc{W}= \{ w_{22} \in L(\St^ {\perpi})^{[s]}: w_{22}= y^{\#}y+ z, \mbox{ where } z=z^ {\#} \geq 0 \}.$$ 
\end{prop}

\begin{dem} Take $J=\matriz{J_1}{0}{0}{J_2}  \ms{\St}{\St^{\perpi}}$ a signature operator for $\HH$ as in \eqref{signature}. 
	
	Then, $JW=\matriz{J_1w_{11}}{J_1w_{12}}{J_2w_{12}^{\#}}{*} \ms{\St}{\St^{\perp}}.$ Then, by \cite[Theorem 2.1]{Baidiuk}, there exists a completion $W$ of $W^0$ if and only if $R(J_1w_{12}) \subseteq R(\vert J_1 w_{11}\vert^ {1/2})$ or equivalently, proceeding as in Theorem \ref{ThmShortedKrein}, $R(w_{12}) \subseteq R(d).$
	
	In this case, by \cite[Theorem 2.1]{Baidiuk}, any selfadjoint operator completion of $JW^0$ admits the representation 
	$$JW=\matriz{J_1w_{11}}{J_1w_{12}}{J_2w_{12}^{\#}}{J_2w_{22}} \ms{\St}{\St^{\perp}},$$ with $J_2w_{22}=f^*uf+z$ and $f$ the reduced solution of $J_1w_{12}=\vert J_1w_{11} \vert^{1/2} x,$ $u$  the partial isometry corresponding to the polar decomposition of $J_1w_{11},$ and $z \in L(\St^{\perp})^+.$ Then, as in the proof of Theorem \ref{ThmShortedKrein}, if $y$ is the unique bounded solution of the equation $w_{12}=dx,$ the set of completions of  $W^0$ has the form $W=\matriz{w_{11}}{w_{12}}{w_{12}^{\#}}{w_{22}}$ with $w_{22}= y^{\#}y+ z, \mbox{ where } z=z^ {\#} \geq 0.$
\end{dem}

\begin{obs} Any completion $W$ of  an incomplete block operator $W^0$ of the form \eqref{WmatrixIncomp} has the same $\St$-compression: $W_{[\St]}.$ Moreover, for any completion $W,$ $W_{[\St]} \leq W.$ 
\end{obs}

\subsection{Comparison with other notions of Schur complement in Krein spaces}

In \cite{XavierMary}, Mary proved that any weakly regular operator $B\in L(\KK,\HH)$ (i.e., any operator such that $\ol{R(B)}$ and $N(B)$ are regular subspaces) admits a (unique) closed Moore-Penrose inverse. That is, there exists an operator $B^{\dagger}: Dom(B^{\dagger})=R(B) \ [\dotplus] \ R(B)^{\perpi} \subseteq \HH \ra \KK$ such that $BB^{\dagger}$ is a symmetric projection from $Dom(B^{\dagger})$ onto $R(B)$ with nullspace $N(B)$ and $B^{\dagger}B$ is a symmetric projection from $\KK$ onto $R(B^{\dagger})=N(B)^{\perpi}$ with nullspace $N(B)$ (see \cite[Corollary 2.9 and Lemma 2.2]{XavierMary}).

Let $W \in L(\HH)^{[s]}$ and let $\St$ be a regular subspace of $\HH$. Suppose that $W$ is represented as in \eqref{WmatrixI} and $w_{11}=dd^{\#}$ is a polar factorization of $w_{11}.$ If $\ol{R(w_{11})}$ is regular, then  $d, d^{\#}$ and $w_{11}$ are weakly regular, therefore, there exist $d^{\dagger}, (d^{\#})^{\dagger}$ and $w_{11}^{\dagger},$ which are weakly regular and $(d^{\#})^{\dagger}=(d^{\dagger})^{\#}$ (see \cite[Theorem 2.8 and Theorem 2.15]{XavierMary}).

Suppose that $W$ is $\St$-weakly complementable. Then $R(w_{12}) \subseteq R(d)$ (see Theorem \ref{ThmShortedKrein}) and $\!$ $d(d^{\dagger}w_{12})=w_{12}.$ Therefore, $d^{\dagger}w_{12} \in L(\St^{\perpi}, \KK)$ is the unique solution of the equation $dx=w_{12}.$ Thus, by Theorem \ref{ThmShortedKrein},

$$W_{/ [\St]}=\matriz{0}{0}{0}{w_{22}-(d^{\dagger}w_{12})^{\#}d^{\dagger}w_{12}} \ms{\St}{\St^{\perpi}}.$$

If in addition $R(w_{11})$ is closed, then $R(d)$ is regular. Thus, $(d^{\dagger})^{\#} d^{\dagger}w_{12}$ is well defined, $w_{11}^{\dagger}=(d^{\#})^{\dagger} d^{\dagger} \in L(\St)$ and 
$$W_{/ [\St]}=\matriz{0}{0}{0}{w_{22}-w_{12}^{\#} w_{11}^{\dagger} w_{12}}=W_{/ [\St]}^{XM},$$ 
where $W_{/ [\St]}^{XM}$ is the Schur complement of $W$ to $\St$ as defined by Mary. See \cite[Theorem 2.20]{XavierMary}.

For a positive operator $W$ in a Hilbert space $\HH$ and a closed subspace $\St \subseteq \HH,$  Pekarev  \cite{Pekarev} showed that the Schur complement $W_{/ \St}$ of $W$ to $\St$ can be expressed as $W_{/ \St}=W^{1/2}(I - P_\M) W^{1/2}$ where $\M=\ol{W^{1/2}(\St)}$. In \cite{SchurKrein}, Pekarev's result was taken as an inspiration to extend the concept to the more general Krein space setting.  In that paper a (bounded) selfadjoint operator $W$ is said to have the {{unique factorization property}} (UFP) if for any two Bogn\'ar-Kr\'amli factorizations of $W=D_1D_1^{\#}=D_2D_2^{\#},$ there is an isomorphism $U$ such that $D_1=D_2U.$

For $W \in L(\HH)^{[s]}$ with the UFP and a closed subspace $\St \subseteq \HH,$ consider $\M=\ol{D^{\#}(\St)}$ and suppose that $\M$ is a regular subspace of $\KK.$ The {{Schur complement}} of $W$ to $\St$ is then defined in that paper as
$$W^{M-MP}_{/ [\St]}=D(I-Q)D^{\#},$$
where $Q$ is the selfadjoint projection onto $\M.$

The next result shows that, when $W$ is $\St$-complementable the regularity of $\M$ can be omitted. If $P_{\M  {\mathbin{\!/\mkern-3mu/\!}} \M^{\perpi}}$ is the {\emph{projection-like}} operator with domain $\M \ [\dotplus] \ \M^{[\perp]}$, range $\M$ and nullspace $\M^{[\perp]}$, then the operator $D(I-P_{\M  {\mathbin{\!/\mkern-3mu/\!}} \M^{\perpi}})D^{\#}$ is well defined and bounded. In this case, 
$$D(I-P_{\M  {\mathbin{\!/\mkern-3mu/\!}} \M^{\perpi}})D^{\#}=W_{/ [\St]}.$$
Since any bounded selfadjoint operator $W$  can be written in the form $W=DD^{\#}$ with $D: \KK \ra \HH$ injective, $\KK$ a Krein space, it follows that $W$ need not have the UFP.

\begin{prop} \label{Prop312} Let $W \in L(\HH)^{[s]}$ and let $\St$ be a closed subspace of $\HH.$ Suppose that $W=DD^{\#}$ with $D \in L(\KK,\HH),$ $N(D)=\{0\}$ and $W$ is $\St$-complementable. Let $\M=\ol{D^{\#}(\St)}.$ Then
	$$\WS=D(I-P_{\M  {\mathbin{\!/\mkern-3mu/\!}} \M^{\perpi}})D^{\#}.$$
\end{prop}

\begin{dem} Since $W$ is $\St$-complementable, we have that
	$\HH= \St + W^{-1}(\St^{\perpi}).$ Suppose that $W=DD^{\#}$ with $D \in L(\KK,\HH)$ and $N(D)=\{0\}.$ Then
	$$R(D^{\#})=D^{\#}(\St) \ [+] \ R(D^{\#}) \cap D^{\#}(\St)^{\perpi}.$$ Furthermore, the sum is direct because $\{0\}=R(D^{\#})^{\perpi} \supseteq D^{\#}(\St)^{\perpi} \cap  \ol{D^{\#}(\St)}.$ 
	Therefore, 
	$$R(D^{\#})=D^{\#}(\St) \ [\dotplus] \ R(D^{\#}) \cap D^{\#}(\St)^{\perpi} \subseteq \M \ [ \dotplus ] \
	\M^{\perpi} .$$
	Let $T := P_{ \M  {\mathbin{\!/\mkern-3mu/\!}} \M^{\perpi} } D^{\#};$ since $R(D^{\#}) \subseteq Dom(P_{ \M  {\mathbin{\!/\mkern-3mu/\!}} \M^{\perpi} } ) = \M \ [ \dotplus ] \ \M^{\perpi},$ $T$ is well defined. Let $Q$ be any projection onto $\St$ such that $WQ=Q^{\#}W.$ Then, for every $x \in \HH,$
	$$Tx=TQx + T(I-Q)x.$$ Since $Qx \in \St,$ $TQx=P_{ \M  {\mathbin{\!/\mkern-3mu/\!}} \M^{\perpi} } D^{\#}Qx=D^{\#}Qx.$ Also, $T(I-Q)x=0$ because $D^{\#}(I-Q)x \in R(D^{\#}(I-Q))=D^{\#}N(Q) \subseteq D^{\#}(W^{-1}(\St^{\perpi}))=R(D^{\#}) \cap D^{\#}(\St)^{\perpi} \subseteq \M^{\perpi}=N(T).$ 
	Therefore, 
	$$T=P_{ \M  {\mathbin{\!/\mkern-3mu/\!}} \M^{\perpi} } D^{\#}=D^{\#}Q \in L(\HH).$$
	Thus, by Corollary \ref{Cor310W},  $$W_{ / [\St]}=W(I-Q)=DD^{\#}(I-Q)=D(I-P_{\M  {\mathbin{\!/\mkern-3mu/\!}} \M^{\perpi}})D^{\#}.$$
\end{dem}

\section*{Acknowledgements}
We thank Professor M.A.~Dritschel who read the original version of this paper and gave us useful advise. 

We thank the anonymous referee for carefully reading our manuscript and helping us to improve
this paper with several useful comments.

Maximiliano Contino and Alejandra Maestripieri were supported by CONICET PIP 0168.  The work of Stefania Marcantognini was done during her stay at  the Instituto Argentino de Matem\'atica with an appointment funded by the CONICET. She is greatly grateful to the institute for its hospitality and to the CONICET for financing her post.

%\section*{References}

\end{document}